\documentclass[12pt]{amsart}
\usepackage{enumitem, color, amsthm, amssymb, amsmath}
\usepackage{tikz}
\usepackage{verbatim}
\usepackage{placeins}
\usepackage{hyperref}
\usepackage{cleveref}
\usepackage{enumerate}
\usepackage{geometry}
\usepackage{aligned-overset}
\usepackage{fullpage}
\usepackage{multirow}
\usepackage{array}
\usepackage{bookmark}
\usepackage{longtable}
\usepackage{caption}
\usepackage{xcolor}
\usepackage{float}
\usepackage{subcaption}
\usepackage{graphicx}
\newtheorem{theorem}{Theorem}
\newtheorem{lemma}[]{Lemma}
\newtheorem{corollary}[]{Corollary}
\newtheorem{proposition}[]{Proposition}

\theoremstyle{definition}

\newtheorem{example}{Example}

\theoremstyle{remark}
\newtheorem{remark}{Remark}

\numberwithin{equation}{section}

\def \l {\lambda}

\def\C{\mathbb C}

\def\K{\mathbf K}

\def\DM{\textrm{DM}}

\def \HD{\textrm{HD}}

\def \kmr{\textrm{kmr}}
\def \DMn{\textnormal{DM}}

\def\Tr{{\rm Tr}}

\def\Z{{\mathbb Z}}

\def\Q{{\mathbb Q}}

\def\K{\mathbb{K}}

\def\SL{\mathrm{SL}}

\def\N{\mathbb N}

\def\({\left(}
\def\){\right)}
\def \l {\lambda}

\def\HDn{\textnormal{HD}}

\def\C{\mathbb{C}}

\def\Z{\mathbb{Z}}
\def\Q{\mathbb{Q}}

\newcommand{\fp}
{\mathbb{F}_p}

\newcommand{\fq}
{\mathbb{F}_q}

\newcommand{\fphat}
{\widehat{\mathbb{F}_{p}^{\times}}}

\newcommand{\fqhat}
{\widehat{\mathbb{F}_{q}^{\times}}}

\newcommand*\HYPERskip{&}
\catcode`,\active
\newcommand*\pFq{
\begingroup
\catcode`\,\active
\def ,{\HYPERskip}%
\doHyper
}
\catcode`\,12
\def\doHyper#1#2#3#4#5{%
\, _{#1}F_{#2}\left[\begin{matrix}#3 \smallskip \\  #4\end{matrix} \; ; \; #5\right]%
\endgroup
}

\catcode`,\active

\catcode`\,12

\catcode`,\active

\catcode`\,12

\catcode`,\active

\catcode`\,12

\begin{document}

\title{On Some Hypergeometric Modularity Conjectures of Dawsey and McCarthy}

\author{Brian Grove}
\address{Mathematics Department\\
  West Virginia University Institute of Technology\\
  Beckley, West Virginia}
\email{brian.grove@mail.wvu.edu}

\date{}

\begin{abstract}
In recent work, the author, in collaboration with Allen, Long, and Tu, developed the Explicit Hypergeometric Modularity Method (EHMM), which establishes the modularity of a large class of hypergeometric Galois representations in dimensions two and three. One important application of the EHMM is the construction of an explicit family of eta-quotients, which we call the $\K_{2}$ functions, from the hypergeometric background. In this article, we introduce an analogous family of eta-quotients, which we call the $\K_{3}$ functions. These $\K_{3}$ functions are constructed using the theory of weight one cubic theta functions originally developed by Jonathan and Peter Borwein. We then use the $\K_{3}$ functions in the EHMM to resolve several hypergeometric modularity conjectures of Dawsey and McCarthy. Further, we provide applications to special $L$-values of the $\K_{3}$ functions and to the study of generalized Paley graphs.
\end{abstract}

\maketitle 
\section{Introduction}\label{sec:intro}

A central goal in number theory is to establish the modularity of motivic Galois representations, predicted by the Langlands program. For example, the modularity of elliptic curves over $\Q$ \cite{wild3adic,TaylorWiles,Wiles} has inspired many important advances in mathematics, such as the resolution of the Sato--Tate conjecture for non-CM elliptic curves over $\Q$ \cite{unconditionalST,conditionalST3,conditionalST1,conditionalST2}. One important class of motivic Galois representations is hypergeometric Galois representations introduced by Katz \cite{Katz90,Katz09}. The traces of these hypergeometric Galois representations involve certain character sums, such as the $H_{p}$ function defined in \eqref{eq:Hpdef}. These hypergeometric character sums naturally arise when counting points on many algebraic varieties over finite fields, such as the Legendre family of elliptic curves \cite{Win3X,koike2dim,Ono98}
\begin{equation*}
E_{\l}:y^{2} = x(x-1)(x-\l)
\end{equation*}
and the Dwork family of hypersurfaces \cite{BarmanDwork,GoodsonDwork,McCarthy-Dwork,Salerno-Dwork,Yu-Dwork} 
\begin{equation*}
X_{\l}^{d}: x_{1}^{d} + x_{2}^{d} + \cdots + x_{d}^{d} = d \l x_{1}x_{2} \cdots x_{d}.
\end{equation*}
Further, the explicit nature of these hypergeometric representations has led to an abundance of new identities and transformations for special $L$-values of many holomorphic modular forms \cite{hmm1,hmm2,ENRosenK2,ENRosenK1}, such as the $\K_{2}$ functions in \eqref{eq:K2def}.

Another important application of hypergeometric Galois representations lies in graph theory. For example, take a fixed finite field $\fq$, with $q \equiv 1 \pmod{4}$ a prime power, as the vertex set of a graph. Then draw an edge between two vertices $a,b$ if and only if $a-b$ is a nonzero square in $\fq$. Graphs constructed in this way are called Paley graphs \cite{jonesPaley,gPaley,originalPaley}. Counting certain properties of Paley graphs and related objects, such as the number of complete subgraphs of a small order, has been accomplished by several authors \cite{bg1,bg4,bg3,bg2,DMpaley,evansPaley,LPSX,ms1,ms2} using hypergeometric character sums. The desire to match these hypergeometric character sums with coefficients of modular forms has motivated the papers \cite{hmm1,hmm2,Dawsey-McCarthy} and this article. The necessary hypergeometric background is as follows.

Fix a positive integer $n$ and consider the multisets $\alpha = \{a_{1}, \ldots, a_{n}\}$ and $\beta =\{b_{1}=1,b_{2}, \ldots, b_{n}\}$ of rational numbers with $b_{2}, \ldots, b_{n} \notin \Z_{\leq 0}$. The collection $\HD = \{\alpha,\beta\}$ is called a \textit{hypergeometric datum} and we will assume that $\HD$ is \textit{primitive}, so $a_{i}-b_{j} \not\in \Z$ for $1 \leq i, j \leq n$. Further, the least positive common denominator of $\HD$ is defined as $M = \textrm{lcd}(\alpha \cup \beta)$. Now let $p \equiv 1 \pmod{M}$ be a prime, $\fp^{\times}$ the multiplicative group for the finite field $\fp$, $\widehat{\fp^{\times}}$ the associated character group, and $\zeta_{n}$ a fixed primitive $n$-th root of unity.

Now for each prime ideal $\mathfrak{p}$ coprime to $M$ we can associate to the residue field $\kappa_{\mathfrak{p}}:= \Z[\zeta_{M}]/\mathfrak{p}$ of size $q$ a character via the $M$-th residue symbol. Then the traces of hypergeometric Galois representations, discussed in \Cref{thm:Katz}, involve the hypergeometric character sums

\begin{equation}\label{eq:Hpdef}
H_{p}\left[\begin{matrix} a_{1} & a_{2} & \ldots & a_{n}\smallskip \\  &b_{2}&\ldots&b_{n} \end{matrix} \; ; \; \l \; ; \, \mathfrak{p} \right] := \frac{1}{1-p} \sum_{\chi \in \widehat{\kappa_{\mathfrak{p}}^{\times}}} \chi((-1)^{n}\lambda) \prod_{i=1}^{n} \frac{g(\iota_{\mathfrak{p}}(a_{j})\chi)g(\iota_{\mathfrak{p}}(-b_{j})\overline{\chi})}{g(\iota_{\mathfrak{p}}(a_{j}))g(\iota_{\mathfrak{p}}(-b_{j}))},
\end{equation}
where $\l \in \kappa_{\mathfrak{p}}^{\times}$,

\begin{equation}\label{eq:iotadef}
\iota_{\mathfrak{p}}\left(\frac{i}{M}\right)(x):= \left(\frac{x}{\mathfrak{p}}\right)_{M}^{i} \equiv x^{(q-1)\frac{i}{M}} \pmod{\mathfrak{p}}
\end{equation} for all $x \in \Z[\zeta_{M}]$, and
$$g(\chi) = \sum_{x \in \fp^{\times}} \chi(x) \zeta_{p}^{x}$$ is the Gauss sum associated to the multiplicative character $\chi \in \fqhat$. Further, we denote the Legendre symbol as $\chi_{d} = (d/\cdot)$ throughout. The $H_{p}$ function was originally introduced by McCarthy \cite{McCarthy} and further developed by Beukers, Cohen, Mellit \cite{bcm}, and others. 

\begin{remark}
The shorthand notation $H_{p}(\HD;\l;\mathfrak{p})$ will be used for the $H_{p}$ function in \eqref{eq:Hpdef} when the hypergeometric datum $\HD$ is clear.
\end{remark}

This article concerns several new hypergeometric modularity results, meaning identities between the $H_{p}(\HD;\l;\mathfrak{p})$ function and the coefficients $a_{p}(f_{\HD}^{\sharp})$ of a Hecke eigenform $f_{\HD}^{\sharp}$ at primes $p \equiv 1 \pmod{M}$. Throughout, $f_{\HD}$ will denote a cusp form on a congruence subgroup of $\SL_{2}(\Z)$ computed from the datum $\HD$ and $f_{\HD}^{\sharp}$ will denote a Hecke eigenform constructed from the cusp form $f_{\HD}$. In particular, we study a list of 15 hypergeometric modularity conjectures made by Dawsey and McCarthy in Table 2 of \cite{Dawsey-McCarthy} for hypergeometric data of the form 
\begin{equation}\label{eq:DMdata}
\HD_{\DM}(u,v) = \{\{1/u,1/v,(v-1)/v\},\{1,1,1\}\}
\end{equation}
at $\l = 1$. In particular, these conjectures take the general form
\begin{equation}\label{eq:DMform}
H_{p}(\HD_{\DM}(u,v);1;\mathfrak{p}) = \psi_{(u,v)}(\mathfrak{p}) \cdot a_{p}\left(f_{\HD_{\DM}(u,v)}^{\sharp}\right)
\end{equation}

for primes $p \equiv 1 \pmod{M}$, where $M$ is the least common denominator of $1/u$ and $1/v$, $\psi_{(u,v)}(\mathfrak{p})$ is an appropriate character on the Galois group $G(M):= \textrm{Gal}(\overline{\Q(\zeta_{M})}/\Q(\zeta_{M}))$, and $f_{\HD_{\DM}(u,v)}^{\sharp}$ is a Hecke eigenform which depends on the pair $(u,v)$. The precise $(u,v)$ pairs considered by Dawsey and McCarthy are 
$$(3,2),(6,2), (8,2), (3,3), (4,3), (6,3), (3,4), (4,4),$$
$$(6,4), (5,5), (3,6), (4,6), (6,6), (2,10), (2,12).$$
Now consider the three $v = 3$ cases, so the pairs $(u,v) \in \{(3,3), (4,3), (6,3)\}$. The main result resolves the conjectures of the form \eqref{eq:DMform} in Table 2 of \cite{Dawsey-McCarthy} for these three pairs with $v = 3$ and two additional cases corresponding to the $(u,v)$ pairs $(2,3)$ and $(12,3)$.

\begin{theorem}\label{thm:K3cases}
Consider the hypergeometric data $\HDn_{\DMn}(u,v)$ from \eqref{eq:DMdata} and let $f_{\HDn_{\DMn}(u,v)}^{\sharp}$ be an explicit Hecke eigenform of weight three constructed from the EHMM. If $(u,v) \in \{(2,3), (3,3), (4,3), (6,3), (12,3)\}$ then
$$H_{p}(\HDn_{\DMn}(u,v);1;\mathfrak{p}) = \psi_{(u,v)}(\mathfrak{p}) \cdot a_{p}\left(f_{\HDn_{\DMn}(u,v)}^{\sharp}\right)$$
at every prime ideal $\mathfrak{p}$ in $\Z[\zeta_{M}]$ above each fixed prime $p \equiv 1 \pmod{M}$, where $M$ is the least common denominator of $1/u$ and $1/v$, for the $\psi_{(u,v)}(\mathfrak{p})$ and $f_{\HDn_{\DMn}(u,v)}^{\sharp}$ values in the table below.
\begin{table}[ht]
\begin{center}
    \begin{tabular}{|c|c|c|c|}
    \hline
    \textrm{$(u,v)$ pair} & $M$ & $\psi_{(u,v)}(\mathfrak{p})$ & $f_{\HDn_{\DMn}(u,v)}^{\sharp}$  \\ \hline
$(2,3)$ & $6$ & $1$ & $f_{12.3.c.a}$\\
\hline
$(3,3)$ & $3$ & $1$ & $f_{27.3.b.b}$\\ 
\hline
$(4,3)$ & $12$ & $\iota_{\mathfrak{p}}(1/4)(-3)$ & $f_{16.3.c.a}$\\
\hline
$(6,3)$ & $6$ & $1$ & $f_{108.3.c.b}$\\
\hline
$(12,3)$ & $12$ & $\iota_{\mathfrak{p}}(1/4)(-3)$ & $f_{432.3.g.e}$\\
\hline
    \end{tabular}
\end{center}
\end{table}

Recall the character $\iota_{\mathfrak{p}}(a)(x)$ is defined in \eqref{eq:iotadef}, and note the $L$-functions and modular forms database (LMFDB) labels are used for all modular forms.

\end{theorem}

\begin{remark}
Note that $\iota_{\mathfrak{p}}(1/4)(-3) \cdot a_{p}(f_{16.3.c.a}) = a_{p}(f_{36.3.d.a})$ at primes $p \equiv 1 \pmod{12}$ which matches the conjecture for the $(u,v) = (4,3)$ pair in Table 2 of \cite{Dawsey-McCarthy}.
\end{remark}

An interesting application of \Cref{thm:K3cases} to counting Paley graphs is as follows. A natural modification of a Paley graph is to replace quadratic residues with $k$-th power residues, called generalized Paley graphs \cite{gPaley}. Fix an integer $k \geq 2$. Now let $q$ be a prime power such that $q \equiv 1 \pmod{k}$ if $q$ is even or $q \equiv 1 \pmod{2k}$ if $q$ is odd. Then define $S_{k}$ to be the order $\frac{q-1}{k}$ subgroup of $\fq^{\times}$ containing the $k$-th power residues. The generalized Paley graph of order $k$, denoted by $G_{k}(q)$, is the graph with vertex set $\fq$, where $ab$ is an edge if and only if $a - b \in S_{k}$. Further, let $\mathcal{K}_{4}(G_{3}(p))$ denote the number of complete subgraphs of order four in $G_{3}(p)$. Now, combining the $(u,v) = (3,3)$ case of \Cref{thm:K3cases} with the discussions in Section 8 of \cite{DMpaley} establishes the following result, conjectured by Dawsey and McCarthy.

\begin{corollary}
If $p \equiv 1 \pmod{6}$ then 
\begin{equation*}
\begin{split}
\mathcal{K}_{4}(G_{3}(p)) = 
 & \, \frac{p(p-1)}{2^{3} \cdot 3^{7}} \bigg[p^{2} +(85- 5p) a_{p}(f_{27.2.a.a})+ 10a_{p}(f_{27.2.a.a})^2+12   a_{p}(f_{27.3.b.b})\\&-55p+316 \bigg].
\end{split}
\end{equation*}
\end{corollary}

Now, many hypergeometric modularity results in the literature rely on various tricks to pinpoint the Hecke eigenform $f_{\HD}^{\sharp}$. For example, several results \cite{Ahlgren-Ono-CalabiYau, lennon1, McCarthy-Papanikolas, Mor3F2} have used techniques such as the Eichler--Selberg trace formula, various identities for the $H_{p}$ function, properties of Gr\"{o}ssencharacters associated to CM modular forms, and known connections between the $H_{p}$ function and elliptic curves.
However, these techniques typically involve fixing a carefully chosen hypergeometric datum and varying the parameter in the associated geometric object, such as an elliptic curve or K3 surface. See \cite{Dawsey-McCarthy} for a more comprehensive survey of many known hypergeometric modularity results.

The lack of a general method for establishing hypergeometric modularity results led the author to establish the Explicit Hypergeometric Modularity Method (EHMM) in collaboration with Allen, Long, and Tu \cite{hmm1,hmm2}. An important feature of the EHMM is an explicit formula for the target modular form $f_{\HD}^{\sharp}$ in many cases when the datum $\HD$ has length three or four with parameter $\l = 1$. See \Cref{sec:EHMM} for more details on the EHMM.

In this article, we construct two new families of weight three eta-quotients with the EHMM, the $\K_{3}$ and $\K_{3}^{\kmr}$ functions. These functions are defined as
\begin{equation}\label{eq:K3def1}
\K_{3}(r,1)(\tau):= \eta(\tau)^{9-12r}\eta(3 \tau)^{12r-3}
\end{equation}
for $r$ in $\mathbb{S}_{3} := \left\{i/12 \, | \, 1 \leq i \leq 11 \right\}$ and 
\begin{equation}\label{eq:K3kummerdef}
\K_{3}^{\rm{kmr}}(r)(\tau):= \eta(\tau)^{1-12r}\eta(3 \tau)^{12r+5},
\end{equation}
for $r$ in $\mathbb{S}_{3}^{\kmr} := \{i/12 \, | \, i \geq 1\}$,
where 
\begin{equation}\label{eq:etadef}
\eta(\tau) := q^{1/24} \prod_{n=1}^{\infty} (1-q^{n})
\end{equation}
is the Dedekind-eta function for $q = e^{2 \pi i \tau}$, where $\tau$ lies in the complex upper-half plane $\mathcal{H} := \{\tau \in \C \, | \, \textrm{Im}(\tau) > 0\}$.

\begin{remark}
The \textrm{kmr} exponent indicates the eta-quotients in \eqref{eq:K3kummerdef} are motivated by the Kummer transformations \eqref{eq:k3kummer}.
\end{remark}

Now consider $r = a/b$ in reduced form and define $N_{\K_{3}}(r) = b$. Then the $\K_{3}(r,1)(N_{\K_{3}}(r)\tau)$ functions are holomorphic cusp forms of weight three defined on certain congruence subgroups $\Gamma_{0}(L)$ with precise levels and Nebentypus given in \Cref{K3table}. Further, the completions of the $\K_{3}$ functions to Hecke eigenforms, which are explicit character twists of the $f_{\HD_\DM(u,3)}^{\sharp}$ functions in \Cref{thm:K3cases}, are also provided in the table below.

Further, if $r \in \{1/12,1/6,1/4\}$ the $\K_{3}^{\kmr}$ functions are $\K_{3}$ functions \eqref{eq:k3tok3kmr}, so holomorphic cusp forms. However, the $\K_{3}^{\kmr}$ functions are not holomorphic, in general. These two families of eta-quotients are constructed with the EHMM using the theory of weight one cubic theta functions, first studied by Jonathan and Peter Borwein \cite{Borweincubic}. The completions of the $\K_{3}$ functions are the eigenforms $f_{\HD_\DM(u,3)}^{\sharp}$ in \Cref{thm:K3cases}, up to an explicit character twist. See \Cref{sec:K3} and \Cref{subsec:K3transform} for the constructions of the $\K_{3}$ and $\K_{3}^{\kmr}$ functions, respectively.

\renewcommand{\arraystretch}{1.4}
\begin{table}[ht]
\caption{The Families of $\K_{3}$ Functions}\label{K3table}
\vspace{2mm}
\begin{center}
    \begin{tabular}{|c|c|c|c|c|c|}
    \hline
    \textrm{Families} & \textrm{$r = j/12$ \, values} & $N_{\K_{3}}(r)$ & \textrm{Eigenform(s)} & \textrm{Level} & \textrm{Character} \\ \hline
$(1)$ & $1/2$ & $2$ & $f_{12.3.c.a} = \K_{3}\left(\frac{1}{2},1\right)(2 \tau)$ & $12$ & $\chi_{-3}$ \\ \hline
$(2)$ &$1/3,2/3$ & $3$ & $f_{27.3.b.b} = \K_{3}\left(\frac{1}{3},1\right)(3 \tau)$ & $27$ & $\chi_{-3}$\\
&&&$+ 3 \sqrt{-1} \K_{3}\left(\frac{2}{3},1\right)(3 \tau)$ && \\ \hline
$(3)$ & $1/4,3/4$ & $4$ & $f_{16.3.c.a} = \K_{3}\left(\frac{1}{4},1\right)(4 \tau)$ & $16$ & $\chi_{-1}$ \\ 
&&&$f_{16.3.c.a}(3 \tau) = \K_{3}\left(\frac{3}{4},1\right)(4 \tau)$ && \\ \hline
$(4)$ & $1/6,5/6$ & $6$ & $f_{108.3.c.b} = \K_{3}\left(\frac{1}{6},1\right)(6 \tau)$ & $108$&$\chi_{-3}$  \\ 
&&&$+9\sqrt{-1}\K_{3}\left(\frac{5}{6},1\right)(6 \tau)$ &&\\ \hline
$(5)$ & $1/12,5/12$  & $12$ & $f_{432.3.g.e} = \K_{3}\left(\frac{1}{12},1\right)(12 \tau)$& $432$ & $\chi_{-1}$\\
& $7/12,11/12$& & $+ 3 \sqrt{5} \K_{3}\left(\frac{5}{12},1\right)(12 \tau)$ &&\\ 
&&&$- 3\sqrt{-15}\K_{3}\left(\frac{7}{12},1\right)(12 \tau)$ &&\\ 
&&&$-9\sqrt{-3}\K_{3}\left(\frac{11}{12},1\right)(12 \tau)$ &&\\ \hline
    \end{tabular}
\end{center}
\end{table}
\renewcommand{\arraystretch}{1}
\newpage

The explicit construction of the $\K_{3}$ functions naturally leads to several other applications, such as new formulas for the special $L$-values of the Hecke eigenforms $f_{\HD_\DM(u,3)}^{\sharp}$ at $s=1$. In particular, the following result is shown.

\begin{corollary}\label{cor:LK3short}
Consider the hypergeometric data $\HDn_\DMn(u,3)$, as in \eqref{eq:DMdata}. Now let $f_{\HDn_\DMn(u,3)}^{\sharp}$ be the Hecke eigenform associated to each datum in \Cref{K3table} and $L(f_{\HDn_\DMn(u,3)},1)$ be the associated special $L$-value at $1$, as in \eqref{eq:Lvaluedef}. Then for each $f_{\HDn_\DMn(u,3)}^{\sharp}$ in \Cref{K3table} there exist algebraic numbers $\alpha_{\DMn(u,3)}$ and  $\beta_{\DM(j)} \in \C$ such that 

\begin{equation*}\label{eq:LK3shortformula}
\alpha_{\DMn(u,3)} \cdot \pi \cdot L(f_{\HDn_\DMn(u,3)}^{\sharp},1) = \sum_{j \in (\Z / u\Z)^{\times}} \beta_{\DMn(j)} \cdot \, _{3}P_{2}(\HD_{\K_{3}}(j/u,1);1)
\end{equation*}

where $_{3}P_{2}(\HD_{\K_{3}}(j/u,1);1)$ is the normalized hypergeometric function defined in \eqref{eq:3P2forK3}.

\end{corollary}

\begin{remark}
See \Cref{subsec:LK3} for the explicit version of \Cref{cor:LK3short}.
\end{remark}

 In addition, applying the Kummer transformation \eqref{eq:k3kummer} for the classical hypergeometric functions connected to the $\K_{3}$ functions, similar to the transformation \eqref{eq:3F2Kummer} studied by Rosen \cite{ENRosenK2} in the $\K_{2}$ case, leads to the $\K_{3}^{\kmr}$ functions \eqref{eq:K3kummerdef}, a second new family of weight three eta-quotients. The $\K_{3}^{\kmr}$ functions exhibit a complementary family of modular forms coming from the EHMM when the original hypergeometric data is not necessarily closed under a Kummer transformation, unlike the $\K_{2}$ case.  

The structure for the remainder of this paper is as follows. In \Cref{sec:EHMM}, the EHMM is recalled with examples. Then in \Cref{sec:K3} the $\K_{3}$ functions are constructed, in analogy with the $\K_{2}$ functions. In \Cref{sec:charsums}, several relevant properties of hypergeometric character sums and the associated Galois representations are discussed in preparation for applying the EHMM. Then in \Cref{sec:K3proofs} the main result, \Cref{thm:K3cases}, is proved. In the final section, \Cref{sec:K3apps}, several applications of the $\K_{3}$ functions, such as special $L$-values and the construction of the $\K_{3}^{\kmr}$ functions, are discussed.
\section*{Acknowledgements}
I thank Ling Long, Esme Rosen, and Fang--Ting Tu for helpful comments and discussions. Further, I thank the anonymous referee for many useful suggestions that greatly improved the quality of the exposition.
\section{The Explicit Hypergeometric Modularity Method (EHMM)}\label{sec:EHMM}

We now recall the EHMM, as developed in \cite{hmm1,hmm2}, in preparation for the proof of \Cref{thm:K3cases}. The EHMM includes two main parts, computing the target modular form using hypergeometric information and establishing several congruences between various types of hypergeometric functions.

An illustrative example of the EHMM involves a family of modular forms called the $\K_{2}$ functions, as defined in \cite{hmm1} and further studied in \cite{hmm2,ENRosenK2}. We will now recall the EHMM for this particular family of modular forms. 

\subsection*{Part I: Computing the Modular Form} The first part is an algorithm to go from a hypergeometric datum $\HD$ to a corresponding modular form $f_{\HD}$. Computing the modular form $f_{\HD}$ involves combining the Euler integral formula for certain classical $_{3}F_{2}$ and $_{4}F_{3}$ hypergeometric functions \cite[Equation (3.2)]{Win3X} with the theory of arithmetic triangle groups, in particular the Schwarz map \cite[Section 3.2.2]{Win3X}. This approach to computing modular forms using hypergeometric information and triangle groups was motivated by earlier work of Li, Long, and Tu \cite{LLT,LLT2}.

Consider hypergeometric data of the form
\begin{equation}\label{eq:K2data}
\HD_{\K_{2}}(r,s) = \{\{1/2,1/2,r\},\{1,1,s\}\},
\end{equation}
for $r,s \in \Q$. We now restrict the $(r,s)$ pairs to a particular subset of $\Q^{2}$ in order to ensure convergence in the Euler integral formula, exclude degenerate cases, and guarantee holomorphic congruence cusp forms. The chosen subset is defined as 
\begin{equation}\label{eq:S2def}
\mathbb{S}_{2} = \{(r,s) \in \Q^{2} \, | \, 0 < r < s < 3/2, r \neq 1, s \neq 1/2, 24s \in \Z, 8(r+s) \in \Z \},
\end{equation}
as in \cite{hmm1,hmm2}. The cardinality of $\mathbb{S}_{2}$ is 193. Now if $(r,s) \in \mathbb{S}_{2}$ the well-known Euler integral formula \cite[Equation (3.2)]{Win3X} applied to the $_{3}F_{2}(\HD_{\K_2}(r,s);1)$ functions gives
\begin{equation}\label{eq:3F2integralK2}
_{3}F_{2}(\HD_{\K_{2}}(r,s);1) = \frac{1}{B(r,s-r)} \int_{0}^{1} t^{r}(1-t)^{s-r-1} \cdot \, _{2}F_{1}(\{1/2,1/2\},\{1,1\};t) \, \frac{dt}{t},
\end{equation}
where $B(x,y)$ is the beta function and $_{2}F_{1}(\{1/2,1/2\},\{1,1\};t)$ is a $_{2}F_{1}$ hypergeometric function over the complex numbers with parameter $t$ \cite[Equation (3.1)]{Win3X}. The idea now is to substitute an appropriate Hauptmodul for $t$ in the differential of \eqref{eq:3F2integralK2}. In particular, the Schwarz map \cite[Section 3.2.2]{Win3X}, gives an association between the $_{2}F_{1}(\{1/2,1/2\},\{1,1\};t)$ function and the congruence subgroup $\Gamma(2)$. Therefore, we are led to the modular lambda function 
$$\lambda(\tau) = 16 \frac{\eta(\tau/2)^{8}\eta(2 \tau)^{16}}{\eta(\tau)^{24}},$$
 It is well-known that $\lambda(\tau)$ is a Hauptmodul for $\Gamma(2)$. Now substituting $t = \lambda(\tau)$ in the differential of \eqref{eq:3F2integralK2} leads to 
\begin{equation}\label{eq:K2eval}
\lambda(\tau)^{r}(1-\lambda(\tau))^{s-r-1} \cdot \, \pFq{2}{1}{\frac{1}{2},\frac{1}{2}}{&1}{\lambda(\tau)} q \frac{d \lambda(\tau)}{\lambda(\tau) dq} = 2^{4r-1} \cdot \K_{2}(r,s)(\tau),
\end{equation}
where the $\K_{2}$ functions are defined as 
\begin{equation}\label{eq:K2def}
\K_{2}(r,s)(\tau) := \frac{\eta(\tau/2)^{16s-8r-12}\eta(2 \tau)^{8s+8r-12}}{\eta(\tau)^{24s-30}},
\end{equation}
for $(r,s) \in \mathbb{S}_{2}$. An important property of the $\K_{2}$ functions, established in \cite{hmm1}, is that when $(r,s) \in \mathbb{S}_{2}$ the $\K_{2}(r,s)(N_{\K_2}(r)\tau)$ functions are weight three holomorphic cusp forms of level $N_{\K_2}(r)N_{\K_2}(s-r)$ with quadratic character $(-2^{24s}/\cdot)$, where 
\begin{equation}\label{eq:Nk2def}
N_{\K_2}(r) = \frac{48}{\textrm{gcd}(24r,24)}.
\end{equation}
\begin{remark}
In addition, the $\K_{2}$ functions can be completed to Hecke eigenforms in certain cases, in analogy with the $\K_{3}$ functions in \Cref{K3table}. See Appendix I of \cite{ENRosenK2} for the complete list of such eigenforms.
\end{remark}

The $\K_{2}$ functions have many applications, such as new formulas for special $L$-values of many eigenforms, as recalled in \Cref{subsec:Lk2}. In addition, the $\K_{2}(r,s)$ functions are related to the Dawsey--McCarthy conjectures for the $(u,v)$ pairs $(3,2), (6,2),$ and $(8,2)$. In the $(3,2)$ and $(6,2)$ cases the exponents on the $\K_{2}$ functions are not all integral so the $\K_{2}(1/3,1)(\tau)$ and $\K_{2}(1/6,1)(\tau)$ functions are noncongruence modular forms of weight three. We hope to address these cases in future work, as the EHMM currently relies on holomorphic modular forms defined on congruence subgroups of $\SL_{2}(\Z)$. Further, the result in the $(8,2)$ case is 
\begin{equation}\label{eq:K21/8}
H_{p}(\HD_{\DM}(8,2);1;\mathfrak{p}) = \chi_{16.e}(p) \cdot a_{p}(f_{256.3.c.g}) = a_{p}(f_{64.3.d.a})
\end{equation}
at prime ideals $\mathfrak{p}$ above primes $p \equiv 1 \pmod{8}$, where $\chi_{16.e}(p)$ is a character in the Dirichlet character orbit 16.e in the LMFDB. The result in \eqref{eq:K21/8} is established in Proposition 6.1 of \cite{hmm1} after accounting for the normalizing factor $\mathcal{J}(\HD;\mathfrak{p})$ in equation $(2.5)$ of \cite{hmm1}, as recalled in \eqref{eq:K21/8}.
\begin{remark}\label{rem:K2structure}
The proof of \eqref{eq:K2eval} relies heavily on the close connection between the modular lambda function $\lambda(\tau)$ and the following three weight $1/2$ Jacobi theta functions:
$$\theta_{2}(\tau) = \sum_{n \in \Z} q^{(2n+1)^{2}/8}, \quad \theta_{3}(\tau) = \sum_{n \in \Z} q^{n^{2}/2}, \quad \theta_{4}(\tau) = \sum_{n \in \Z} (-1)^{n} q^{n^{2}/2}.$$
See \eqref{eq:lambda2theta}, \eqref{eq:keythetaidentity}, and Section 2.4 of \cite{LLT} for the key connections between these theta functions and $\lambda(\tau)$.
\end{remark}
\subsection*{Part II: Several congruences modulo \texorpdfstring{$p^2$}{p squared}}
The other main part of the EHMM consists of combining various congruences modulo $p^2$, for primes $p \equiv 1 \pmod{M}$. Assume that $p$ is an ordinary prime, meaning the quantities $a_{p}(f_{\HD}^{\sharp})$ and $ H_{p}(HD;1;\mathfrak{p})$ do not vanish modulo $p$. The first key observation is that proving an equality between the $H_{p}(\HD;1;\mathfrak{p})$ function and the Fourier coefficients $a_{p}(f_{\HD}^{\sharp})$, up to an appropriate character twist, is equivalent to proving a congruence modulo $p^2$. This reduction is possible due to the Deligne bound on $a_{p}(f_{\HD}^{\sharp})$ and a Weil-type bound on the $H_{p}(\HD;1;\mathfrak{p})$ function due to Katz \cite{Katz90}. 

The EHMM relies on three main congruences. The first two congruences are bridged by the truncated hypergeometric series 
$$F(\HD;1)_{p-1} = \pFq{n}{n-1}{a_{1}&a_{2}&\cdots&a_{n}}{&b_{2}&\cdots&b_{n}}{1}_{p-1} = \sum_{k=0}^{p-1} \frac{\prod_{j=1}^{n}(a_{j})_{k}}{k! \prod_{j=2}^{n}(b_{j})_{k}}$$
when $n = 3,4$, where $(a)_{0} = 1$ and $(a)_{k} := a(a+1) \cdots (a+k-1)$ for $k > 0$. In particular, the first congruence relates the $H_{p}(\HD;1;\mathfrak{p})$ function to the $F(\HD;1)_{p-1}$ function by the Gross-Koblitz formula and various tools from $p$-adic analysis mainly coming from \cite{LR}. Then the second congruence relates the $F(\HD;1)_{p-1}$ function to a $p$-adic unit root $\mu_{p}(\HD;1)$ by the residue sum technique, as used in \cite{Allen,LTYZ}. Now the final congruence relates the $p$-adic unit root $\mu_{p}(\HD;1)$ to the coefficients $a_{p}(f_{\HD}^{\sharp})$, up to an appropriate character, modulo $p^2$ using the modularity of hypergeometric Galois representations and the explicit computation of $f_{\HD}^{\sharp}$ from Part I of the EHMM. See \cite{hmm1} for more details on the EHMM.

We now recall a special case of the EHMM. See Theorem 2.1 of \cite{hmm1} for the EHMM in full generality. In particular, we restrict to the case of $n = 3$ with $\HD = \{\{a_{1},a_{2},a_{3}\},\{1,1,b_{3}\}\}$ such that
$$\gamma(\HD):= 1 + b_{3} - \sum_{i=1}^{3} a_{i} \neq 1.$$
The restrictions above guarantee that the hypergeometric Galois representations discussed in \Cref{thm:EHMMrestricted} include the cases with hypergeometric data $\HD_{\DM}(u,v)$. Now the EHMM has two main assumptions. The first assumption involves the existence of a holomorphic cusp form $f_{\HD}$ computed from the Euler integral formula, such as the $\K_{2}$ functions in \eqref{eq:K2eval} or the $\K_{3}$ functions in \eqref{eq:K3def}. In particular, this cusp form must also be completed to a Hecke eigenform $f_{\HD}^{\sharp}$, as in the $\K_{3}(1/12,1)(12\tau)$ case from \Cref{ex:TpforK3example}.

Then the second assumption involves checking if the Galois representation $\rho_{\{\HD;1\}}$ attached to the datum $\HD$, as defined in \Cref{thm:Katz}, is extendable from the Galois group $G(M)$ to the absolute Galois group $G_{\Q}:= \textrm{Gal}(\overline{\Q}/\Q)$. Luckily, the traces of hypergeometric Galois representations involve hypergeometric character sums, such as the $H_{p}$ function in \eqref{eq:Hpdef}. Therefore, the second condition can be phrased at the level of hypergeometric character sums.

Further, a strength of the EHMM is that every step is explicit. For example, an explicit formula for $f_{\HD}$ involving the differential of the Euler integral formula and exact formulas for various character twists built from the hypergeometric data and the leading coefficient of the chosen Hauptmodul $t(\tau)$ are provided. The specialization of the EHMM under the restrictions discussed above is as follows.
\begin{theorem}[\cite{hmm1}]\label{thm:EHMMrestricted}
Consider $\alpha^{\flat} = \{a_{1},a_{2}\}$, where $0 < a_{1} \leq a_{2} < 1$ and $\beta^{\flat} = \{1,1\}$ with $a_{3},b_{3}$ such that $0 < a_{3} < b_{3} \leq 1$ and $a_{2} < b_{3}$. Let $\HDn = \{\alpha^{\flat} \cup \{a_{3}\}, \beta^{\flat} \cup \{b_{3}\}\}$ be primitive and $M = \rm{lcd}(\HDn)$. Suppose that the following two conditions are satisfied.

$(1)$ There exists a modular function $t = C_{1}q + O(q^{2}) \in \Z [[q]]$ with $C_{1} \neq 0$ such that 
$$f_{\HDn}(q):= C_{1}^{-a_{3}} \cdot t(q)^{a_{3}}(1-t(q))^{b_{3}-a_{3}-1} \cdot \, _{2}F_{1}(\alpha^{\flat}, \beta^{\flat}; t(q)) q \frac{d t(q)}{t(q) dq}$$
is a congruence holomorphic cuspform of weight three. Further, for each prime $p \equiv 1 \pmod{M}$ there exist integers $\tilde{b}_{p}$ such that $T_{p}(f_{\HDn}) = \tilde{b}_{p} \cdot f_{\HDn}$, where $T_{p}$ is the $p$-th Hecke operator.

$(2)$ For any prime ideal $\mathfrak{p} \in \Z[\zeta_{M}]$ above a fixed prime $p \equiv 1 \pmod{M}$,
$$-\iota_{\mathfrak{p}}(a_{3})(C_{1})^{-1} J(\iota_{\mathfrak{p}}(a_{3}),\iota_{\mathfrak{p}}(b_{3}-a_{3})) \cdot H_{p}(\HDn;1;\mathfrak{p}) \in \Z,$$
where
\begin{equation*}
J(\iota_{\mathfrak{p}}(a),\iota_{\mathfrak{p}}(b)):= \sum_{x \in \kappa_{\mathfrak{p}}} \iota_{\mathfrak{p}}(a)(x) \cdot \iota_{\mathfrak{p}}(b)(1-x)
\end{equation*}
is a Jacobi sum, defined for $n \geq 2$. 

Then there exists a normalized Hecke eigenform $f_{\HDn}^{\sharp}$ built from $f_{\HDn}$, not necessarily unique, such that for each prime $p > 29$ with $p \equiv 1 \pmod{M}$, $\tilde{b}_{p} = a_{p}(f_{\HDn}^{\sharp})$. More explicitly, 
$$a_{p}(f_{\HDn}^{\sharp}) = -\iota_{\mathfrak{p}}(a_{3})(C_{1})^{-1} \cdot J(\iota_{\mathfrak{p}}(a_{3}),\iota_{\mathfrak{p}}(b_{3}-a_{3})) \cdot H_{p}(\HDn;1;\mathfrak{p}).$$
In terms of Galois representations,
$$\rho_{f_{\HDn}^{\sharp}}|_{G(M)} \cong \chi_{\HDn} \otimes \rho_{\{\HDn;1\}},$$
where $\rho_{f_{\HDn}^{\sharp}}$ is the Galois representation associated to the eigenform $f_{\HDn}^{\sharp}$ by Deligne, $\rho_{\{\HDn;1\}}$ is the hypergeometric Galois representation associated to $\HDn$ from \Cref{thm:Katz} by Katz, and 
$$\chi_{\HDn}(\mathfrak{p}):= \iota_{\mathfrak{p}}(a_{3})(C_{1})^{-1} \cdot \prod_{i=1}^{2} \iota_{\mathfrak{p}}(a_{i})(-1).$$
\end{theorem}
%
%
%
\begin{remark}
The EHMM applies to primes $p \equiv 1 \pmod{M}$, in general. The small primes can be checked numerically. In particular, the $p > 29$ bound in \Cref{thm:EHMMrestricted} is written as a general lower bound that can be reduced in many individual cases. 
\end{remark}
We are now ready to discuss an important class of eta-quotients constructed in the EHMM, the $\K_{3}$ functions, which can be viewed as a cubic analog of the $\K_{2}$ functions in \eqref{eq:K2def}. These $\K_{3}$ functions are used to build the Hecke eigenforms $f_{\HD_{\DM}(u,v)}^{\sharp}$ in \Cref{thm:K3cases}.
\section{The \texorpdfstring{$\K_{3}$}{K3} Functions}\label{sec:K3}

\subsection{Constructing the \texorpdfstring{$\K_{3}$}{K3} Functions}\label{subsec:K3construction}

The computation of the $\K_{2}$ functions in \eqref{eq:K2eval} relied heavily on the relationships between the Hauptmodul $\lambda(\tau)$ and the three weight $1/2$ Jacobi theta functions $\theta_{2}(\tau), \theta_{3}(\tau),$ and $\theta_{4}(\tau)$, as mentioned in \Cref{rem:K2structure}. For example, the modular lambda function $\lambda(\tau)$ satisfies 
\begin{equation}\label{eq:lambda2theta}
\lambda(\tau) = \left(\frac{\theta_{2}(\tau)}{\theta_{3}(\tau)} \right)^{4}.
\end{equation}
Further, the three-term identity 
\begin{equation}\label{eq:keythetaidentity}
\theta_{3}(\tau)^{4} = \theta_{2}(\tau)^{4} + \theta_{4}(\tau)^{4},  
\end{equation}
originally established by Jacobi, is used in constructing the $\K_{2}$ functions and in many other applications of theta functions, such as for the theories of elliptic integrals and arithmetic-geometric means \cite{aar,BB}. See \cite{BB, Zagier-modularform} for more details on the Jacobi theta functions given in \eqref{eq:keythetaidentity}. Now, an interesting question is if there are other Hauptmoduln $t(\tau)$ constructed from the Schwarz map and a corresponding family of theta functions such that equations analogous to \eqref{eq:lambda2theta} and \eqref{eq:keythetaidentity} hold.

Such relations were predicted for certain weight one cubic theta functions by the work of Ramanujan \cite{Ramanujan-pi} and later established by the Borwein brothers and Garvan \cite{BBG-Ramanujan, Borweincubic, BBG}. We begin our discussion of the cubic case with hypergeometric functions and the Schwarz map. Consider the hypergeometric data 
\begin{equation}\label{eq:HDk3}
\HD_{\K_{3}}(r,1):= \{\{1/3,2/3,r\},\{1,1,1\}\}
\end{equation}
for $(r,1) \in \mathbb{S}_{3}$, where 
\begin{equation}\label{eq:S3def}
\begin{split}
\mathbb{S}_{3}&:= \{(r,s) \in \Q^{2} \, | \, 0 < r < 1, 12r \in \Z, s = 1\}\\
&= \left\{\left(\frac{j}{12},1\right) \, \bigg| \, 1 \leq j \leq 11 \right\}.
\end{split}
\end{equation}
The choice of $s=1$ is made to ensure the hypergeometric data corresponds to a modular form. Further, recall the group $\Gamma_{0}(3)$ has an elliptic point of order three, unlike the group $\Gamma(2)$. This fact provides some explanation for why the set $\mathbb{S}_{3}$ is more restrictive, so smaller, compared to the $\mathbb{S}_{2}$ set \eqref{eq:S2def}.
\begin{remark}
The hypergeometric data $\HD_{\K_{3}}(r,1)$ can be viewed as a cubic analog of \eqref{eq:K2data}. Note the $s = 1$ assumption in \eqref{eq:S3def} is more restrictive then the assumptions on $s$ in the $\mathbb{S}_{2}$ set \eqref{eq:S2def} associated to the $\K_{2}$ functions. Further, the condition $0 < r < 1$ is required for convergence in the Euler integral formula, and the $12r \in \Z$ condition ensures holomorphic modular forms.
\end{remark}
Now if $(r,1) \in \mathbb{S}_{3}$ the Euler integral formula applied to the $_{3}F_{2}(\HD_{\K_{3}}(r,1);1)$ functions gives
\begin{equation}\label{eq:3F2integralK3}
_{3}F_{2}(\HD_{\K_{3}}(r,1);1) = \frac{1}{B(r,1-r)} \int_{0}^{1} t^{r}(1-t)^{-r} \cdot \, _{2}F_{1}(\{1/3,2/3\},\{1,1\};t) \frac{dt}{t},
\end{equation}
in analogy with \eqref{eq:3F2integralK2}. Then the Schwarz map, as discussed in Section 3.2.2 of \cite{Win3X}, gives an association between the $_{2}F_{1}(\{1/3,2/3\},\{1,1\};t)$ function and the congruence subgroup $\Gamma_{0}(3)$. In this case, a Hauptmodul for $\Gamma_{0}(3)$ is the modular function
\begin{equation}\label{eq:t3def}
t_{3}(\tau) = 27 \frac{\eta(3 \tau)^{9}}{\left(3 \eta(3 \tau)^{3} + \eta\left(\frac{\tau}{3}\right)^{3}\right)^{3}},
\end{equation}
as discussed in \cite{hmm1, Borweincubic,BBG-Ramanujan, BBG}. Now the analogs of the weight $1/2$ Jacobi theta functions $\theta_{2}(\tau), \theta_{3}(\tau),$ and $\theta_{4}(\tau)$ in the cubic setting are the weight one cubic theta functions $a(\tau), b(\tau),$ and $c(\tau)$ \cite{BBG} given as follows. Let $\zeta_{3} = e^{2 \pi i/3}$. Then define  $a(\tau),b(\tau),$ and $c(\tau)$ as
\begin{equation}\label{eq:abcdef}
\begin{split}
a(\tau) := \sum_{(n,m) \in \Z^{2}} q^{n^{2}+nm+m^{2}} = \frac{3 \eta(3 \tau)^{3}+\eta(\tau/3)^{3}}{\eta(\tau)},\\
\\
b(\tau):= \sum_{(n,m) \in \Z^{2}} \zeta_{3}^{m-n}q^{n^{2}+nm+m^{2}} = \frac{\eta(\tau)^{3}}{\eta(3 \tau)},\\
\\
c(\tau):= \sum_{(n,m) \in \Z^{2}} q^{(n+1/3)^{2}+(n+1/3)(m+1/3)+(m+1/3)^{2}} = 3 \frac{\eta(3 \tau)^{3}}{\eta(\tau)}.
\end{split}
\end{equation}
The explicit representations of $a(\tau), b(\tau),$ and $c(\tau)$ as eta-quotients in \eqref{eq:abcdef} lead to the identity 
\begin{equation}\label{eq:t3withabc}
t_{3}(\tau) = \left(\frac{c(\tau)}{a(\tau)}\right)^{3},
\end{equation}
in analogy with \eqref{eq:lambda2theta}. Further, the work of the Borwein brothers \cite{Borweincubic} and their later joint paper with Garvan \cite{BBG} provides two different proofs of the three-term identity
\begin{equation}\label{eq:abckeyidentity}
a(\tau)^{3} = b(\tau)^{3} + c(\tau)^{3},
\end{equation}
which can be viewed as a cubic analog of \eqref{eq:keythetaidentity}. We now require a few more results before substituting $t=t_{3}(\tau)$ in the differential of \eqref{eq:3F2integralK3}, in analogy with \eqref{eq:K2eval}. First, combining the identities \eqref{eq:t3withabc} and \eqref{eq:abckeyidentity} gives 
\begin{equation}\label{eq:1minust3}
1 - t_{3}(\tau) = \left(\frac{b(\tau)}{a(\tau)} \right)^{3}. 
\end{equation}
The remaining two necessary results involve evaluating the $_{2}F_{1}(\{1/3,2/3\},\{1,1\};t)$ function and the logarithmic derivative $dt/t$ in the differential of \eqref{eq:3F2integralK3}, both at $t = t_{3}(\tau)$. First, the work of the Borwein brothers in \cite{Borweincubic} gives
\begin{equation}\label{eq:2F1att3}
_{2}F_{1}(\{1/3,2/3\},\{1,1\};t_{3}(\tau)) = a(\tau).
\end{equation}
Now, the logarithmic derivative evaluation follows from the lemma below.
\begin{lemma}\label{lem:t3der}
Let $t_{3}(\tau), a(\tau),b(\tau),$ and $c(\tau)$ be as defined in \eqref{eq:t3def} and \eqref{eq:abcdef}. Further, define the differential operator $D(f) := (2 \pi i)^{-1} \cdot df/d\tau$. Then
$$D(t_{3}(\tau)) = \frac{c^{3}(\tau)b^{3}(\tau)}{a^4(\tau)}$$ so
\begin{equation}\label{eq:t3logder}
 D(\log(t_{3}(\tau))) = \frac{b^{3}(\tau)}{a(\tau)}.
\end{equation}
\end{lemma}

\begin{proof}
It is shown that 
\begin{equation}\label{eq:derofS}
D\left(\frac{c(\tau)}{a(\tau)}\right) = \frac{1}{3} \cdot \frac{c(\tau)b^{3}(\tau)}{a^{2}(\tau)}
\end{equation}
in Section 2.4 of \cite{LLT}. The results then follow by \eqref{eq:t3withabc} and \eqref{eq:derofS}.
\end{proof}

Now substituting $t = t_{3}(\tau)$ in the differential of \eqref{eq:3F2integralK3} using \eqref{eq:t3def}, \eqref{eq:abcdef}, \eqref{eq:t3withabc}, \eqref{eq:abckeyidentity}, \eqref{eq:1minust3}, and \eqref{eq:2F1att3} gives 
\begin{equation}\label{eq:K3eval}
\begin{split}
&t_{3}(\tau)^{r}(1-t_{3}(\tau))^{-r} \cdot \, _{2}F_{1}(\{1/3,2/3\},\{1,1\};t_{3}(\tau)) \cdot D( \log(t_{3}(\tau)))\\
&= b^{3(1-r)} \cdot c^{3r}\\
&= 27^{r} \K_{3}(r,1)(\tau)
\end{split}
\end{equation}
where 
\begin{equation}\label{eq:K3def}
\K_{3}(r,1)(\tau):= \eta(\tau)^{9-12r}\eta(3 \tau)^{12r-3}.
\end{equation}
Now, to show that the $\K_{3}(r,1)(\tau)$ functions are holomorphic cuspforms of weight three for $(r,1) \in \mathbb{S}_{3}$, we rely on a well-known result. The following result checks if an eta-quotient, which means a function built from quotients of the Dedekind-eta function \eqref{eq:etadef}, is a holomorphic modular form on a congruence group $\Gamma_{0}(L)$ with explicit Nebentypus.
\begin{theorem}\rm(\cite{Ono-WebofModularity})\label{thm:etacheck}
    Suppose there is some $L > 4$ such that $f(\tau) =  \displaystyle\prod_{d \, | \, L} \eta(d \tau)^{r(d)}$ with $r(d)$ integers. If
    \begin{itemize}
    \item[$(1)$] $k:= \frac{1}{2}  \displaystyle\sum_{d \, | \, L} r(d) \in \Z,$

    \item[$(2)$] $\displaystyle\sum_{d \, | \, L} d \cdot r(d) \equiv \sum_{d \, | \, L} \frac{L}{d} \cdot r(d) \equiv 0 \pmod{24},$ and
   
    \item[$(3)$] $\frac{L}{24} \displaystyle\sum_{d \, | \, L} \frac{\textrm{gcd}(c,d)^{2} \cdot r(d)}{\textrm{gcd}(c, \frac{L}{c}) \cdot cd}$ is nonnegative (resp. positive) for each divisor $c \, | \, L$
    \end{itemize}
    
    then $f(\tau) \in M_{k}(\Gamma_{0}(L), \chi)$ (resp. $f(\tau) \in S_{k}(\Gamma_{0}(L), \chi)$), where $\chi(\cdot) = \bigg(\frac{(-1)^{k}\prod_{d \, | \, L} d^{\, r(d)}}{\cdot} \bigg)$ is a quadratic character.  
\end{theorem}
Condition $(1)$ of \Cref{thm:etacheck} shows that the $\K_{3}(r,1)(\tau)$ functions are of weight three. Then the explicit level $L(r)$, change of variables $\tau \mapsto N_{\K_{3}}(r)\tau$, and character $\chi = \chi(r)$ are determined by minimizing $d$ and $L$ in conditions $(2)$ and $(3)$ of \Cref{thm:etacheck}. For example, the change of variables $\tau \mapsto N_{\K_{3}}(r)\tau$ is computed as follows.

Define the quantity $N_{\K_3}(r)$ as the smallest positive integer such that the $\K_{3}(r,1)(N_{\K_{3}}(r)\tau)$ function in \eqref{eq:K3def} satisfies condition $(2)$ of \Cref{thm:etacheck}. Then $N_{\K_{3}}(r)$ satisfies $24r \cdot N_{\K_{3}}(r) \equiv 0 \pmod{24}$. Equivalently, we desire the smallest positive integer $N_{\K_3}(r)$ such that $N_{\K_3}(r) \cdot r$ is a positive integer. If $r = a/b$ is in reduced form, then taking $N_{\K_3}(r) = b$ gives the desired minimization. Similar arguments involving conditions $(2)$ and $(3)$ of \Cref{thm:etacheck} provide the minimal levels and the precise characters for all eleven $\K_{3}(r,1)(N_{\K_{3}}(\tau))$ functions in \Cref{K3table}. 

We now discuss the cases in which the $\K_{2}$ and $\K_{3}$ functions can be completed to Hecke eigenforms in detail.

\subsection{From the \texorpdfstring{$\K_{3}$}{K3} Functions to Hecke Eigenforms}\label{subsec:galoisK2andK3}

The cases in which the cusp forms constructed from the EHMM can be explicitly completed to eigenforms are studied in \cite{hmm1,hmm2,ENRosenK2} for the $\K_{2}$ functions. We now explain a similar construction for the $\K_{3}$ functions in \Cref{K3table}. In particular, we start by constructing the eigenform $f_{\HD_\DM(u,3)}^{\sharp}$ from \Cref{thm:K3cases} with the $\K_{3}$ functions in the case $u = 12$.

\begin{example}\label{ex:TpforK3example}
The eigenform $f_{\HD_\DM(12,3)}^{\sharp}$ is constructed as a particular combination of the $\K_{3}$ functions in the last row of \Cref{K3table}. In particular, the work of Rosen \cite{ENRosenK2} implies that the constants in this combination of $\K_{3}$ functions can be computed using the relevant $T_{p}$ Hecke operators on the cuspforms in family $(5)$ of \Cref{K3table}. Observe that \Cref{K3table} suggests the subspace of cuspforms $S_{3}(\Gamma_{0}(432), \chi_{-1})$ generated by the $f_{j} = \K_{3}(j/12,1)(12 \tau) \in q^{j}(1 + \Z[[q^{12}]])$ functions for $j=1,5,7,11$ as a starting point. We now compute the $T_{p}$ Hecke operators for $p = 2,3,5,7,11$ to construct the eigenform $f_{432.3.g.e}$ from the $f_{j}$ functions. It is straightforward to check that $T_{2}(f_{j}) = T_{3}(f_{j}) = 0$ for $j=1,5,7,11$ in this case. The nontrivial actions of the Hecke operators $T_{5}, T_{7},$ and $T_{11}$ are recorded in the table below.
$$
  \begin{array} {|c||c|c|c|c|}
    \hline
    &f_1&f_5&f_7&f_{11}      \\\hline 
    T_5&45f_{5}&f_{1}&9f_{11}&5f_{7}  \\\hline
    T_7&-135f_{7}&-27f_{11}&f_{1}&5f_{5}  \\\hline
    T_{11}&-243f_{11}&-27f_{7}&9f_{5}&f_{1}  \\\hline
  \end{array}
$$
The above table showcases many relations between the $T_{p}$ operators for $p \in \{5,7,11\}$, such as 
\begin{equation}\label{eq:HeckerelationsK3ex}
T_{5}^{2} = 45, \quad  T_{7}^{2} = -135, \quad T_{11}^{2} = -243, \quad \textrm{and} \quad T_{5}T_{7} = 5T_{11}.
\end{equation}
The EHMM then predicts a normalized Hecke eigenform $f_{\HD}^{\sharp}(\tau)$ of the form
$$f_{\HD}^{\sharp}(\tau) = \sum_{j \in (\Z/12 \Z)^{\times}} c_{j} \cdot \K_{3}(j/12,1)(12 \tau),$$
for constants $c_{1},c_{5},c_{7},c_{11}$. The normalized assumption gives $c_{1} = 1$. Now the constants $c_{5},c_{7},$ and $c_{11}$ are computed from the Hecke operators $T_{5}, T_{7},$ and $T_{11}$, respectively. In particular, the constants $c_{5}, c_{7},$ and $c_{11}$ satisfy the relations $c_{5}^{2} = 45, c_{7}^{2} = -135, c_{11}^{2} = -243,$ and $c_{5}c_{7} = 5 c_{11}$ by \eqref{eq:HeckerelationsK3ex}. Note the eigenform $f_{\HD}^{\sharp}(\tau)$ is not unique as the choices of $c_{5}$ and $c_{7}$ uniquely determine $c_{11}$. For example, choosing $c_{5} = 3 \sqrt{5}$ and $c_{7} = -3 \sqrt{-15}$ gives $c_{11} = -9 \sqrt{-3}$, which agrees with \Cref{K3table}. Further, the eigenform $f_{\HD}^{\sharp}(\tau)$ for this choice is labeled as $f_{432.3.g.e}$ in the LMFDB.
\end{example}

Observe how the eigenform completion was constructed by computing the relevant Hecke operators in the above example for the $\K_{3}(1/12,1)(12 \tau)$ case. A similar computation gives all the eigenform completions for the $\K_{3}$ functions in \Cref{K3table}. In particular, the constants in the combinations of $\K_{3}$ functions for families $(2)$ and $(4)$ follow from the following Hecke operator computations.
\begin{table}[ht]
\caption{The action of $T_{2}$ on the subspace of $S_{3}(\Gamma_{0}(27), \chi_{-3})$ generated by the $g_{j} = \K_{3}(j/3,1)(3 \tau)$ functions for $j=1,2$.}\label{family2table}
\begin{center}
    \begin{tabular}{|c|c|c|}
    \hline
    &$g_{1}$&$g_{2}$\\
    \hline
    $T_{2}$&$-9g_{2}$&$g_{1}$\\
    \hline
    \end{tabular}
\end{center}
\end{table}
\begin{table}[ht]
\caption{The action of $T_{5}$ on the subspace of $S_{3}(\Gamma_{0}(108), \chi_{-3})$ generated by the $h_{j} = \K_{3}(j/6,1)(6 \tau)$ functions for $j=1,5$.}\label{family4table}
\begin{center}
    \begin{tabular}{|c|c|c|}
    \hline
    &$h_{1}$&$h_{5}$\\
    \hline
    $T_{5}$&$-81h_{5}$&$h_{1}$\\
    \hline
    \end{tabular}
\end{center}
\end{table}

 We now return to hypergeometric character sums and condition $(2)$ of \Cref{thm:EHMMrestricted} in preparation for the proof of \Cref{thm:K3cases} in \Cref{sec:K3proofs}.

\section{Hypergeometric Character Sums and Extendability}\label{sec:charsums}

The traces of hypergeometric Galois representations involve the $H_{p}$ function, as mentioned earlier. However, it is sometimes more convenient to use the following normalization of the $H_{p}$ function
\begin{equation}\label{eq:Pdef}
\mathbb{P}\left[\begin{matrix} a_{1} & a_{2} & \ldots & a_{n}\smallskip \\  &b_{2}&\ldots&b_{n} \end{matrix} \; ; \; \l \; ; \, \mathfrak{p} \right] := \prod_{i=2}^{n} J(\iota_{\mathfrak{p}}(a_{i}),\iota_{\mathfrak{p}}(b_{i}-a_{i})) \cdot \, H_{p}\left[\begin{matrix} a_{1} & a_{2} & \ldots & a_{n}\smallskip \\  &b_{2}&\ldots&b_{n} \end{matrix} \; ; \; \l \; ; \, \mathfrak{p} \right].
\end{equation}
 The $\mathbb{P}$ function in \eqref{eq:Pdef} is inspired by the hypergeometric character sums of Greene \cite{Greene} and Katz \cite{KatzESDE} and was originally defined as the $_{n}\mathbb{P}_{n-1}$ function in \cite{Win3X}. Further, the shorthand notation $\mathbb{P}(\HD;\l;\mathfrak{p})$ will be used for \eqref{eq:Pdef}, in analogy with the $H_{p}$ case. The traces of hypergeometric Galois representations involve a particular character twist of the $\mathbb{P}(\HD;1;\mathfrak{p})$ function. In particular, the precise foundational results of Katz on hypergeometric Galois representations are recalled below.
\begin{theorem}[Katz \cite{Katz90, Katz09}]\label{thm:Katz}  Let  $\ell$ be a prime. Given a primitive hypergeometric datum $\HDn = \{\{r_1,\cdots,r_n\},\{q_{1}=1,q_2,\cdots,q_n\}\}$ with $M = \textnormal{lcd}(\HDn)$ and 
$$\gamma(\HDn) = -1 + \sum_{i=1}^{n} (q_{i}-r_{i}),$$
for
any $\l \in \Z[\zeta_M,1/M]\smallsetminus \{0\}$ the following hold. 
\begin{itemize}
\item [i).]There exists an $\ell$-adic Galois representation $\rho_{\{\HDn;\l\}}: G(M)\rightarrow GL(W_{\l})$ unramified almost everywhere such that at each nonzero
prime ideal $\mathfrak{p}$ of  $\Z[\zeta_M,1/(M\ell \l)]$ of norm $N(\mathfrak{p})=|\Z[\zeta_M)]/\mathfrak{p}|$ 
\begin{equation*}\label{eq:Tr1} \Tr \rho_{\{\HDn;\l\}}(\textnormal{Frob}_\mathfrak{p})= (-1)^{n-1}  \iota_{\mathfrak{p}}(r_1) (-1)\cdot \mathbb{P}(\HDn; 1/\l;\mathfrak{p}),  
\end{equation*} 
where $\textnormal{Frob}_\mathfrak{p}$ denotes the  geometric  Frobenius conjugacy class
of $G(M)$ at $\mathfrak{p}$. 

\item[ii).] When $\l\neq 1$,  the dimension $d := dim_{\overline \Q_\ell}W_{\l}$ equals $n$ and all roots of the characteristic polynomial of $\rho_{\{\HDn;\l\}}(\textnormal{Frob}_\mathfrak{p})$  are algebraic numbers and have the same absolute value $N(\mathfrak{p})^{(n-1)/2}$ under all archimedean embeddings. If $\rho_{\{\HDn;\l\}}$ is self-dual, namely isomorphic to any of its complex conjugates, then $W_{\l}$ admits a non-degenerate alternating (resp. symmetric)   bilinear pairing if
$n$ is even and $\gamma(\HDn)\in\Z$ (resp. otherwise).
\item[iii).] When $\l=1$, in the self-dual case the dimension is $n-1$. All roots of the Frobenius eigenvalues at $\wp$  have absolute value less or equal to $N(\mathfrak{p})^{(n-1)/2}$.  If $\rho_{\{\HDn;\l\}}$ is self-dual, it contains a subrepresentation that admits a non-degenerate alternating (resp. symmetric)   bilinear pairing if
$n$ is even and $\gamma(\HDn)\in\Z$ (resp. otherwise). For this subrepresentation, the roots of the characteristic polynomial of $\mathrm{Frob}_{\mathfrak{p}}$ have absolute value exactly $N(\mathfrak{p})^{(n-1)/2}$.
\end{itemize}
\end{theorem}
Now since the traces of Katz representations involve the $H_{p}$ function in many cases condition $(2)$ of \Cref{thm:EHMMrestricted} is equivalent to establishing certain character sum identities between the $H_{p}(\HD_{\K_{3}}(r,1);1;\mathfrak{p})$ functions for all Galois conjugates of $r$ in these cases. This equivalence is possible by the following well-known result for proving the extendability of Galois representations.
\begin{proposition}[\cite{SerreAbelian}] \label{prop:extendableprop}
Let $M$ be a positive integer. Assume $\rho$ is a semi-simple finite dimensional $\ell$-adic representation of $G(M)$ which is isomorphic to $\rho^{\tau}$ for each $\tau \in \textrm{G}_{\Q}$, then $\rho$ is extendable to $\textrm{G}_{\Q}$. Equivalently, for each nonzero prime ideal $\mathfrak{p}$ of $\Z[\zeta_{M}]$ unramified for $\rho$, we have $\Tr \, \rho(\textnormal{Frob}_{\mathfrak{p}}) \in \Z$.
\end{proposition}
In light of \Cref{thm:Katz} and \Cref{prop:extendableprop}, the traces of the $\rho_{\{\HD;1\}}$ representations for Galois conjugates are related to each other by transformation formulas for the $H_{p}$ function. There are many character sum identities for the $H_{p}(\HD_{\K_{3}}(r,1);1;\mathfrak{p})$ functions, many of which were originally established by Greene \cite{Greene}. For example, see Section 3.2 of \cite{DMpaley} for a list of many useful identities. The identity we use for the $H_{p}$ function is
\begin{equation}\label{eq:Hptransform}
H_{p}\left[\begin{matrix} a&1-a&c\smallskip \\  &1&1\end{matrix} \; ; \; 1 \; ; \, \mathfrak{p} \right] = \iota_{\mathfrak{p}}(a)(-1) \cdot H_{p}\left[\begin{matrix} a&1-a&1-c\smallskip \\  &1&1\end{matrix} \; ; \; 1 \; ; \, \mathfrak{p} \right],
\end{equation}
at primes $p \equiv 1 \pmod{M}$ where $M = \textnormal{lcd}(a,b,c)$ and $a,b,c \in \Q \setminus \Z$.
\begin{remark}
The identity \eqref{eq:Hptransform} is derived as follows. First, taking $A,B,C \in \fphat$ to be nontrivial, $B = \overline{A}$, and $D = E = \varepsilon$ the trivial character in $(3.18)$ of \cite{DMpaley} gives the result for the finite field $_{3}F_{2}(1)$ functions of Greene. Then the identity follows by the conversion between the $H_{p}$ function and Greene's finite field hypergeometric function in Proposition 2.5 of \cite{McCarthy}.
\end{remark}
We are interested specifically in the hypergeometric data $\HD_{\K_{3}}(r,1)$ here so taking $a = 1/3, b = 2/3,$ and $c = r$ in \eqref{eq:Hptransform} gives 
\begin{equation}\label{eq:HptransformK3}
H_{p}\left[\begin{matrix} \frac{1}{3}&\frac{2}{3}&r\smallskip \\  &1&1\end{matrix} \; ; \; 1 \; ; \, \mathfrak{p} \right] = H_{p}\left[\begin{matrix} \frac{1}{3}&\frac{2}{3}&1-r\smallskip \\  &1&1\end{matrix} \; ; \; 1 \; ; \, \mathfrak{p} \right]
\end{equation}
at primes $p \equiv 1 \pmod{\textrm{lcd}(1/3,r)}$ for $r \in \Q \setminus \Z$. Now observe that the identities \eqref{eq:Hptransform} and \eqref{eq:HptransformK3} are involutions. Therefore, \eqref{eq:HptransformK3} suffices to establish condition $(2)$ of \Cref{thm:K3cases} for families $(1)-(4)$ in \Cref{K3table}. 

Further, for family $(5)$ the conjugate pairs $(1/12,11/12)$ and $(5/12,7/12)$ of $r$ values can be directly related by \eqref{eq:HptransformK3}. However, relating $r = 1/12$ to $r = 5/12$ (or $r = 7/12$) cannot be done with \eqref{eq:HptransformK3} directly since such a transformation is not an involution. To our knowledge, the pairs $r = 1/12$ and $r = 5/12$ cannot be related directly with known identities for the $H_{p}$ function with parameter $\l = 1$.

Therefore, for family $(5)$ in \Cref{K3table} we use a contradiction argument, as this orbit has size four.
\begin{proposition}\label{prop:K3galoiscase}
    For every prime $p \equiv 1 \pmod{12}$ and prime ideal $\mathfrak{p} \in \Z[\zeta_{12}]$ above $p$ we have 
    $$\iota_{\mathfrak{p}}(1/4)(-3) \cdot H_{p}(\textnormal{HD}_{\K_{3}}(1/12,1);1;\mathfrak{p}) = a_{p}(f_{432.3.g.e}).$$
\end{proposition}
\begin{proof}
First, observe that $\iota_{\mathfrak{p}}(1/4)(-3)$ is a quadratic character of $G(M)$ under the given assumptions. Further, the $a_{p}(f_{432.3.g.e})$ values are integers with absolute value less than or equal to $2p$, by the Deligne bound. Note the left-side takes values in $\Z[\zeta_{12}]$, as discussed in \cite{bcm,Win3X}, so it can be written in the form 
$$a_{1} + a_{2} \sqrt{-1} + a_{3} \sqrt{3} + a_{4} \sqrt{-3}$$
for $a_{j} \in \frac{1}{2}\Z$. Now there are four unique ways to embed the ring $\Z[\zeta_{12}]$ in $\Q_{p}$. In particular, fix an embedding $\sigma(\zeta_{12}) \in \Q_{p}$. Then all embeddings take the form $\sigma_{j}:\zeta_{12} \mapsto \sigma(\zeta_{12})^{j}$ for $j=1,5,7,11$. Now the supercongruences in \cite{hmm1}, such as Theorem 2.3, imply that 
\begin{equation}\label{eq:congmodp2}
\sigma_{j}(a_{1} - a_{p}(f_{432.3.g.e}) + a_{2} \sqrt{-1} + a_{3} \sqrt{3} + a_{4} \sqrt{-3}) \equiv 0 \pmod{p^2}
\end{equation}
for all $j$. Then \eqref{eq:congmodp2} implies that $p^{2}$ divides each of the half integers $a_{1}-a_{p}(f_{432.3.g.e}), a_{2}, a_{3},$ and $a_{4}$. Now suppose that at least one of the half integers $a_{1}, \ldots, a_{4}$ is nonzero. This implies that any complex norm of
$$\iota_{\mathfrak{p}}(1/4)(-3) \cdot \, H_{p}(\HD_{\K_{3}}(1/12,1);1;\mathfrak{p}) - a_{p}(f_{432.3.g.e})$$
is at least $p^2/2$. However, any complex norm of either $\iota_{\mathfrak{p}}(1/4)(-3) \cdot \, H_{p}(\HD_{\K_{3}}(1/12,1);1;\mathfrak{p})$ or $a_{p}(f_{432.3.g.e})$ is less than or equal to $2p$, by the Deligne bound and \Cref{thm:Katz}. This is impossible when $p \geq 13$, so the difference between both sides is zero.
\end{proof}

We are now ready to complete the proof of \Cref{thm:K3cases}.

\section{The Proof of \texorpdfstring{\Cref{thm:K3cases}}{Theorem 2.1}}\label{sec:K3proofs}

We now prove \Cref{thm:K3cases} using the EHMM in \Cref{thm:EHMMrestricted}.

\begin{proof}[Proof of \Cref{thm:K3cases}]
First, note the hypergeometric data $\HD_{\K_{3}}(r,1)$, defined in \eqref{eq:HDk3}, satisfies the initial assumptions in \Cref{thm:EHMMrestricted} when $(r,1) \in \mathbb{S}_{3}$. Further, the $r=1/12$ case is shown in \Cref{prop:K3galoiscase} so the remaining cases to consider correspond to families $(1)-(4)$. Then condition $(1)$ of \Cref{thm:EHMMrestricted} follows by the construction of the $\K_{3}(r,1)(\N_{\K_{3}}(r)\tau)$ functions in \eqref{eq:K3eval} and as discussed in \Cref{sec:intro}  and \Cref{sec:K3}. Note that for families $(1)$ and $(3)$ in \Cref{K3table}, the conjugate orbit contains a single $\K_{3}$ function, up to a change of variables in family $(3)$, which is a Hecke eigenform. In particular, 
$$f_{12.3.c.a} = \K_{3}(1/2,1)(2 \tau) \quad \textrm{and} \quad f_{16.3.c.a} = \K_{3}(1/4,1)(4 \tau)$$ 
for families $(1)$ and $(3)$, respectively. In the remaining cases, families $(2), (4),$ and $(5)$, the actions of the relevant $T_{p}$ Hecke operators on the subspaces of modular forms spanned by the appropriate Galois families of $\K_{3}(r,1)(N_{\K_{3}}(r)\tau)$ functions are computed explicitly in \Cref{subsec:galoisK2andK3}. Further, the Hecke eigenforms $f_{\HD_{\DM}(u,3)}^{\sharp}$ constructed from the $\K_{3}(r,1)(N_{\K_{3}}(r)\tau)$ functions for the pairs $(r,1) \in \mathbb{S}_{3}$ are listed in \Cref{K3table}.

Now consider condition $(2)$ of \Cref{thm:EHMMrestricted}. The remaining conjugate families, mainly families $(1)-(4)$ in \Cref{K3table}, have at most two conjugates. Then \Cref{thm:EHMMrestricted} says that showing the quantity
\begin{equation*}\label{eq:EHMMcondition2}
\begin{split}
&-\iota_{\mathfrak{p}}(r)(1/27) J(\iota_{\mathfrak{p}}(r), \iota_{\mathfrak{p}}(1-r)) \cdot H_{p}(\HD_{\K_{3}}(r,1);1;\mathfrak{p})\\
&= \iota_{\mathfrak{p}}(r)(-1/27) \cdot H_{p}(\HD_{\K_{3}}(r,1);1;\mathfrak{p})
\end{split}
\end{equation*}
is integer-valued for the remaining $(r,1)$ pairs in $\mathbb{S}_{3}$ and any fixed prime ideal $\mathfrak{p} \in \Z[\zeta_{M}]$ above a prime $p \equiv 1 \pmod{M}$ will establish condition $(2)$. It is straightforward to check that $\iota_{\mathfrak{p}}(r)(-1/27) = 1$ for $r = 1/2,1/3,1/6$ at primes $p \equiv 1 \pmod{6}$ and $\iota_{\mathfrak{p}}(1/4)(-1/27) = \iota_{\mathfrak{p}}(1/12)(-1/27) = \iota_{\mathfrak{p}}(1/4)(-3)$ is a quadratic character at primes $p \equiv 1 \pmod{12}$, since $-3$ is a square at these primes. Now \Cref{prop:extendableprop} and \eqref{eq:HptransformK3} imply that the $H_{p}(\HD_{\K_{3}}(r,1);1;\mathfrak{p})$ functions are integer-valued for $r = 1/2,1/3,1/4,$ and $1/6$.

Therefore, the hypotheses in \Cref{thm:EHMMrestricted} are satisfied for these four cases. Then applying \Cref{thm:EHMMrestricted} to the datum $\HD_{\DM}(u,3) = \HD_{\K_{3}}(1/u,1)$ for $u=2,3,4,6$ gives
$$H_{p}(\HD_{\DMn}(u,3);1) = \psi_{(u,3)}(\mathfrak{p}) \cdot a_{p}\left(f_{\HD_{\DMn}(u,3)}^{\sharp}\right),$$
where $\psi_{(u,3)}(\mathfrak{p}) = \iota_{\mathfrak{p}}(1/u)(-1/27)$. The precise $\psi_{(u,3)}(\mathfrak{p})$ and $f_{\HD_{\DM}(u,3)}^{\sharp}$ values, as discussed above, are listed in the statement of \Cref{thm:K3cases}.
\end{proof}

\begin{remark}
The values of the function $H_{p}(\HD;\l;\mathfrak{p})$ in \eqref{eq:Hpdef} do not depend on the ordering of parameters in the separate sets $\alpha$ and $\beta$ of hypergeometric data. However, the presence of the Jacobi sums in \eqref{eq:Pdef} indicates the values of the $\mathbb{P}(\HD;1;\mathfrak{p})$ function do depend on the ordering or parameters.
\end{remark}

We now provide two applications involving the $\K_{3}$ functions and \Cref{thm:K3cases}. The first application involves special $L$-values for the $\K_{3}$ functions and some related transformations.

\section{Special \texorpdfstring{$L$}{L}-Values for the \texorpdfstring{$\K_{3}$}{K3} Functions and Related Transformations}\label{sec:K3apps}

\subsection{A Review of Special \texorpdfstring{$L$}{L}-Values for the \texorpdfstring{$\K_{2}$}{K2} Functions}\label{subsec:Lk2} The explicit construction of the $\K_{2}$ and $\K_{3}$ functions from the Euler integral formula leads to many new formulas for special values of many modular $L$-series. The study of special $L$-values for the $\K_{2}$ functions was initiated in \cite{hmm1,hmm2,ENRosenK2} and was motivated by the related work \cite{special3vals,momentint}. Before stating the main result which links special $L$-values of the $\K_{2}$ functions with hypergeometric functions we recall a normalization of the $_{3}F_{2}(\HD_{\K_2}(r,s);1)$ functions which allows for cleaner statements. Define 
\begin{equation}\label{eq:classicalPforK2}
_{3}P_{2}(\HD_{\K_{2}}(r,s);1) := \pi \cdot B(r,s-r) \cdot \, _{3}F_{2}(\HD_{\K_{2}}(r,s);1),
\end{equation}
where $B(x,y)$ is the beta function, as before.

\begin{remark}
The $_{3}P_{2}(\HD_{\K_{2}}(r,s);1)$ function is a special case of the hypergeometric period functions defined in equation $(3.6)$ of \cite{Win3X}.
\end{remark}

Let $f(\tau)$ be a modular form of weight $k$. It is well-known \cite{K-Z,LLT,Zagier-top-diff} that the special $L$-values $L(f,1)$ can be computed as
\begin{equation}\label{eq:Lvaluedef}
L(f,1) := -2 \pi i \int_{0}^{i \infty} f(t) dt.
\end{equation}
Now integrating the differentials which define the $\K_{2}$ and $\K_{3}$ functions, as in \eqref{eq:3F2integralK2} and \eqref{eq:3F2integralK3}, naturally leads to identities between the $L(f,1)$ and $_{3}P_{2}$ values, when $f$ is a $\K_{2}$ or $\K_{3}$ function (or a completion thereof). First, consider the $\K_{2}$ case. The key result in this case was shown in \cite{hmm2} and is recalled below.

\begin{lemma}[\cite{hmm2}]\label{lem:K2Lvalues-main}
Let $(r,s) \in \mathbb{S}_{2}$, $N = N_{\K_{2}}(r)$ as in \eqref{eq:Nk2def}, and recall the special $L$-values in \eqref{eq:Lvaluedef}. Then 
\begin{equation*}
_{3}P_{2}(\HDn_{\K_{2}}(r,s);1) = 2^{4r-1} N \pi \cdot L(\K_{2}(r,s)(N\tau),1)
\end{equation*}
\end{lemma}

We now illustrate \Cref{lem:K2Lvalues-main} with an example. 

\begin{example}
Consider the datum $\HD_{\K_{2}}(1/4,3/4) = \{\{1/2,1/2,1/4\},\{1,1,3/4\}\}$. Then, completing the $\K_{2}$ function corresponding to this datum leads to the following eigenform
\begin{equation}\label{eq:K2_1/4_ex}
f_{32.3.c.a}(\tau) = \K_{2}(1/4,3/4)(8 \tau) + 4i \cdot \K_{2}(3/4,5/4)(8 \tau),
\end{equation}
where the choice $\sqrt{-1} = i$ is made. Now combining \eqref{eq:K2_1/4_ex} with \Cref{lem:K2Lvalues-main} gives
\begin{equation}\label{eq:K2_1/4_lvalue}
\pi \cdot L(f_{32.3.c.a},1) = \frac{1}{8} \left[\, _{3}P_{2}(\HD_{\K_{2}}(1/4,3/4);1) + i \cdot \, _{3}P_{2}(\HD_{\K_{2}}(3/4,5/4);1) \right].
\end{equation}
\end{example}

Another interesting application of \Cref{lem:K2Lvalues-main} involves new transformations between the $L(\K_{2}(r,s)(N\tau);1)$ special values. This direction was pursued in \cite{hmm2} and by Rosen \cite{ENRosenK2}. Such transformations are induced by known transformation formulas for the $_{3}F_{2}(\HD_{\K_{2}}(r,s);1)$ functions. For example, the identity
\begin{equation}\label{eq:3F2Kummer}
_{3}F_{2}(\HD_{\K_{2}}(r,s);1) = \frac{\Gamma(s)\Gamma(s-r)}{\Gamma(s-1/2)\Gamma(s-r+1/2)} \cdot \, _{3}F_{2}(\HD_{\K_{2}}(1-r,1/2-r+s);1)
\end{equation}
listed in Corollary 3.3.5 of \cite{aar} was originally shown by Kummer, and it holds when both sides converge, where $\Gamma(x)$ is the Gamma function. Now combining \eqref{eq:3F2Kummer} with \Cref{lem:K2Lvalues-main} leads to interesting transformations, such as
\begin{equation}\label{eq:Lvaluetransform}
L(\K_{2}(1/4,3/4)(8 \tau),1) = 2 \sqrt{2} \cdot L(\K_{2}(3/4,1)(8 \tau),1).
\end{equation}
Further, in many cases the transformations for $L$-values of the $\K_{2}$ functions extend to new transformations for the special $L$-values of the eigenforms which are built from the $\K_{2}$ functions in the orbit, such as the $f_{32.3.c.a}$ case in \eqref{eq:K2_1/4_lvalue}. For example, in the $\K_{2}(1/4,3/4)(8 \tau)$ example discussed above the appearance of the $\K_{2}(3/4,1)(8 \tau)$ function in \eqref{eq:Lvaluetransform} indicates a transformation between the values $L(f_{32.3.c.a},1)$ and $L(f_{64.3.c.b},1)$ since
$$f_{64.3.c.b} = \K_{2}(1/4,1)(8 \tau) + 4i \cdot \K_{2}(3/4,1)(8 \tau).$$
The desired relation in this case is
\begin{equation}\label{eq:rosentransform}
L(f_{32.3.c.a},1) = \frac{\sqrt{2}}{2}i \cdot L(f_{64.3.c.b},1),
\end{equation}
established by Rosen in Section 4.2 of \cite{ENRosenK2}. The relation \eqref{eq:rosentransform} is one example of many relations established by Rosen in \cite{ENRosenK2}.

We now proceed to the special $L$-values for the $\K_{3}$ functions.

\subsection{Special \texorpdfstring{$L$}{L}-Values for the \texorpdfstring{$\K_{3}$}{K3} Functions}\label{subsec:LK3}

The first step is to define the $_{3}P_{2}$ function in this setting, in analogy with \eqref{eq:classicalPforK2}. In particular, define 
\begin{equation}\label{eq:3P2forK3}
_{3}P_{2}(\HD_{\K_{3}}(r,1);1) :=  \frac{2 \sqrt{3}}{3}\pi \cdot  B(r,1-r) \cdot \, _{3}F_{2}(\HD_{\K_{3}}(r,1);1).
\end{equation}
The main result that connects the $_{3}P_{2}(\HD_{\K_{3}}(r,1);1)$ values with the special $L$-values of the $\K_{3}(r,1)(\tau)$ functions at $s = 1$ is below.

\begin{lemma}\label{lemma:K3Lvalues}
Let $(r,1) \in \mathbb{S}_{3}$, $N = N_{\K_{3}}(r)$ as in \Cref{K3table}, and recall the special $L$-values in \eqref{eq:Lvaluedef}. Then
\begin{equation*}
\begin{split}
_{3}P_{2}(\HDn_{\K_{3}}(r,1);1) = 2 \cdot 3^{3r-\frac{1}{2}}N\pi \cdot L(\K_{3}(r,1)(N\tau),1)
\end{split}
\end{equation*}
\end{lemma}

\begin{proof}
Let $t_{3} = t_{3}(\tau)$ be the Hauptmodul defined in \eqref{eq:t3def}. Note that $t_{3}(0) = 1$ and $t_{3}(i \infty) = 0$ in the fundamental domain for $\Gamma_{0}(3)$, as discussed in both Section 3.1 and the Appendix of \cite{hmm1}. Now combining \eqref{eq:3F2integralK3}, \eqref{eq:K3eval}, \eqref{eq:Lvaluedef}, and \eqref{eq:3P2forK3} gives 
\begin{equation*}
\begin{split}
_{3}P_{2}(\HD_{\K_{3}}(r,1);1) &= \frac{2\pi}{\sqrt{3}} \cdot \int_{0}^{1} t_{3}^{r}(1-t_{3})^{-r} \cdot \pFq{2}{1}{\frac{1}{3}&\frac{2}{3}}{&1}{t_{3}} \frac{dt_{3}}{t_{3}}\\
&= 2 \cdot 3^{3r-\frac{1}{2}}\pi \cdot L(\K_{3}(r,1)(\tau),1).
\end{split}
\end{equation*}
The result follows after scaling $\tau \mapsto N \tau$.
\end{proof}

An interesting corollary of \Cref{lemma:K3Lvalues} is a list of the $L$-value results corresponding to the eigenforms built from the $\K_{3}$ families $(1)-(5)$ in \Cref{K3table}. The following result is the explicit version of \Cref{cor:LK3short}.

\begin{corollary}\label{cor:LK3vals} Let $_{3}P_{2}(\HDn_{\K_{3}}(r,1);1)$ be as defined in \eqref{eq:3P2forK3}. Then the following results are true. 
\begin{equation*}
12\pi \cdot L(f_{12.3.c.a},1) = \, _{3}P_{2}(\HDn_{\K_{3}}(1/2,1);1)
\end{equation*}

\begin{equation*}
2 \cdot 3^{\frac{3}{2}}\pi \cdot L(f_{27.3.b.b},1) = \,  _{3}P_{2}(\HDn_{\K_{3}}(1/3,1);1) + \sqrt{-1} \cdot \, _{3}P_{2}(\HDn_{\K_{3}}(2/3,1);1)
\end{equation*}

\begin{equation*}
8 \cdot 3^{\frac{1}{4}}\pi \cdot L(f_{16.3.c.a},1) = \, _{3}P_{2}(\HDn_{\K_{3}}(1/4,1);1) = \, \frac{1}{\sqrt{3}} \cdot \, _{3}P_{2}(\HDn_{\K_{3}}(3/4,1);1)
\end{equation*}

\begin{equation*}
12\pi \cdot L(f_{108.3.c.b},1) = \, _{3}P_{2}(\HDn_{\K_{3}}(1/6,1);1) + \sqrt{-1} \cdot \, _{3}P_{2}(\HDn_{\K_{3}}(5/6,1);1)
\end{equation*}

\begin{equation*}
\begin{split}
8 \cdot 3^{\frac{5}{4}}\pi \cdot L(f_{432.3.g.e},1) &= \sqrt{3} \left[\, _{3}P_{2}(\HDn_{\K_{3}}(1/12,1);1) + \sqrt{5} \cdot \, _{3}P_{2}(\HDn_{\K_{3}}(5/12,1);1) \right]\\
&-\sqrt{-3} \left[\sqrt{5} \cdot \, _{3}P_{2}(\HDn_{\K_{3}}(7/12,1);1) + \, _{3}P_{2}(\HDn_{\K_{3}}(11/12,1);1) \right]
\end{split}
\end{equation*}
\end{corollary}

We now proceed to various transformations and applications for the $\K_{3}$ functions, such as an analog of the Kummer transformation \eqref{eq:3F2Kummer}.

\subsection{Transformations and Applications of the \texorpdfstring{$\K_{3}$}{K3} functions}\label{subsec:K3transform}

Two important sources of symmetry for the $_{3}F_{2}(\HD_{\K_{2}}(r,s);1)$ and $\K_{2}(r,s)$ functions are the Kummer transformation \eqref{eq:3F2Kummer} and the Atkin-Lehner involution $\tau \mapsto -1/\tau$. Recall that the $\K_{2}$ functions are closed under the Kummer transformation \eqref{eq:3F2Kummer}. This closure property for the $\K_{2}$ functions leads to many transformations between the $L(\K_{2}(r,s)(N_{\K_{2}}(r),1)$ values when $(r,s) \in \mathbb{S}_{2}$, as studied in detail by Rosen \cite{ENRosenK2}.

\subsubsection{Kummer Transformations}

It is then natural to ask if an analog of the Kummer transformation \eqref{eq:3F2Kummer} exists for the $_{3}F_{2}(\HD_{\K_{3}}(r,1);1)$ functions when $(r,1) \in \mathbb{S}_{3}$. Now specializing the Kummer transformation in Corollary 3.3.5 of \cite{aar} leads to the following transformations
\begin{equation}\label{eq:k3kummer}
_{3}F_{2}(\HD_{\K_{3}}(r,1);1) = \frac{\Gamma(1-r)}{\Gamma((3-j)/3)\Gamma((3+j)/3-r)} \cdot \, _{3}F_{2}(\HD_{\K_{3},j}(r,1);1)
\end{equation}
for $j=1,2$ fixed, where 
\begin{equation}
\HD_{\K_{3},j}(r,1):= \{\{j/3,j/3,1-r\},\{1,1,(3+j)/3-r\}\}.
\end{equation}
The transformations in \eqref{eq:k3kummer} suggest a bifurcation of the 
datum $\HD_{\K_{3}}(r,1)$ to $\HD_{\K_{3,j}}(r,1)$ for $j=1,2$, unlike the $\K_{2}$ case in \eqref{eq:3F2Kummer}. We now follow the EHMM to compute the modular forms corresponding to the data $\HD_{\K_{3},j}(r,1)$ for $j=1,2$. These corresponding modular forms, which we call the $\K_{3,j}(r,s)$ functions, are constructed below.

First, recall the $t_{3}(\tau)$ and $a,b,c$ functions defined in \eqref{eq:t3def} and \eqref{eq:abcdef}, respectively. Further, let $j \in \{1,2\}$ for the discussion below. Now applying Pfaff transformations for the $_{2}F_{1}$ function \cite[Equation 3.8]{Win3X} and using \eqref{eq:2F1att3} gives
\begin{equation}\label{eq:2F1kummereval}
\begin{split}
\pFq{2}{1}{\frac{j}{3}&\frac{j}{3}}{&1}{\frac{t_{3}}{t_{3}-1}} &= (1-t_{3})^{\frac{j}{3}} \cdot \, \pFq{2}{1}{\frac{1}{3}&\frac{2}{3}}{&1}{t_{3}}\\
&= (1-t_{3})^{\frac{j}{3}} \cdot a.
\end{split}
\end{equation}
Then the Euler integral formula \cite[Equation (3.2)]{Win3X} for this case is 
\begin{equation}\label{eq:Eulerint}
\pFq{3}{2}{\frac{j}{3}&\frac{j}{3}&r}{&1&s}{1} = \frac{1}{B(r,s-r)} \int_{0}^{1} t^{r}(1-t)^{s-r-1} \cdot \, \pFq{2}{1}{\frac{j}{3}&\frac{j}{3}}{&1}{t} \frac{dt}{t},
\end{equation}
defined for $s > r > 0$. Evaluating the integrand of \eqref{eq:Eulerint} at $t = t_{3}/(t_{3}-1)$ leads to
\begin{align*}
&\frac{1}{2 \pi i} \left(t^{r}(1-t)^{s-r-1} \cdot \, \pFq{2}{1}{\frac{j}{3}&\frac{j}{3}}{&1}{t} \frac{dt}{t d \tau} \right)_{t=t_{3}/(t_{3}-1)}\\
&\overset{\eqref{eq:2F1kummereval}}{=} \left(\frac{t_{3}}{t_{3}-1}\right)^{r} \cdot \left(1 - \frac{t_{3}}{t_{3}-1}\right)^{s-r-1} \cdot (1-t_{3})^{\frac{j}{3}} \cdot a \cdot D\left(\log\left(\frac{t_{3}}{t_{3}-1}\right)\right)\\
&\overset{\eqref{eq:t3withabc},\eqref{eq:1minust3}}{=} (-1)^{r} \cdot \left(\frac{c}{b}\right)^{3r-1} \cdot a \cdot \left(\frac{b}{a}\right)^{j+3+3r-3s} \cdot \frac{d}{d\tau} \left(\frac{t_{3}}{1-t_{3}}\right)\\
&\overset{\Cref{lem:t3der}}{=} (-1)^{r} \cdot \left(\frac{c}{b}\right)^{3r-1} \cdot a \cdot \left(\frac{b}{a}\right)^{j+3+3r-3s} \cdot \left(\frac{a}{b}\right)^{6} \cdot \frac{c^{3}b^{3}}{a^{4}}\\
&= (-1)^{r} \cdot a^{3s-3r-j} \cdot b^{j+1-3s} \cdot c^{3r+2}\\
&\overset{\eqref{eq:abcdef}}{=} (-1)^{r} \cdot \left(\frac{3 \eta(3 \tau)^{3}+ \eta(\tau/3)^{3}}{\eta(\tau)}\right)^{3s-3r-j} \cdot \left(\frac{\eta(\tau)^{3}}{\eta(3 \tau)}\right)^{j+1-3s} \cdot \left(3 \frac{\eta(3 \tau)^{3}}{\eta(\tau)} \right)^{3r+2} \\
&= 9 \cdot (-27)^{r} \cdot  \eta(\tau)^{1+4j-12s}\eta(3 \tau)^{(5-j)+9r+3s} \left(3 \eta(3 \tau)^{3} + \eta(\tau/3)^{3}\right)^{3(s-r)-j}.
\end{align*}
Now define 
\begin{equation}\label{eq:K3j}
\K_{3,j}(r,s)(\tau) := \eta(\tau)^{1+4j-12s}\eta(3 \tau)^{(5-j)+9r+3s} \left(3 \eta(3 \tau)^{3} + \eta(\tau/3)^{3}\right)^{3(s-r)-j}.
\end{equation}
Then choosing $s = r + j/3$ in \eqref{eq:K3j} specializes the $\K_{3,j}$ functions to the eta-quotients
\begin{equation}\label{eq:K3kmr}
\K_{3}^{\rm{kmr}}(r)(\tau):= \eta(\tau)^{1-12r}\eta(3 \tau)^{12r+5}.
\end{equation}

\begin{remark}
 It is interesting to note the $\K_{3}^{\textrm{kmr}}(r)(\tau)$ function is independent of $j$. 
\end{remark}

Now \eqref{eq:K3kmr} and the requirement that $r > 0$ in the Euler integral formula suggest the following subset of $\Q$ for which the $\K_{3}^{\textrm{kmr}}(r)(\tau)$ functions are eta-quotients constructed with the EHMM
\begin{equation}\label{eq:S3kmrdef}
\mathbb{S}_{3}^{\kmr} := \{r \in \Q \, | \, r > 0, 12r \in \Z\}.
\end{equation}
It is now natural to ask about connections between the $\K_{3}$ and $\K_{3}^{\kmr}$ functions. First, the relation
\begin{equation*}
\K_{3}(r,1)(\tau) = \left(\frac{\eta(\tau)}{\eta(3 \tau)} \right)^{8} \cdot \K_{3}^{\kmr}(r)(\tau)
\end{equation*}
is clear from \eqref{eq:K3def} and \eqref{eq:K3kmr}. Another interesting question is if any eta-quotients appear as both $\K_{3}$ and $\K_{3}^{\kmr}$ functions. Note that the $\mathbb{S}_{3}$ set, defined in \eqref{eq:S3def}, implies that the smallest exponent on $\eta(\tau)$ in the $\K_{3}(r,1)(\tau)$ function is $-2$ and the largest exponent on $\eta(3 \tau)$ in the $\K_{3}(r,1)(\tau)$ function is $8$. Then applying these bounds to the exponents of the $\K_{3}^{\kmr}(r)(\tau)$ function leads to the bounds $1-12r \geq -2$ and $5+12r \leq 8$ which both give $0 < r \leq 1/4$, since $r$ is positive in $\mathbb{S}_{3}^{\kmr}$. Further, the definition of the $\mathbb{S}_{3}^{\kmr}$ set implies that the $\K_{3}^{\kmr}(\tau)$ function can only be a $\K_{3}(r,1)(\tau)$ function for $r \in \{1/12,1/6, 1/4\}$. Comparing \eqref{eq:K3def} and \eqref{eq:K3kmr} leads to the following three identities
\begin{equation}\label{eq:k3tok3kmr}
\begin{split}
\K_{3}^{\kmr}(1/12)(\tau) &= \K_{3}(3/4,1)(\tau),\\
\K_{3}^{\kmr}(1/6)(\tau) &= \K_{3}(5/6,1)(\tau),  \\ \K_{3}^{\kmr}(1/4)(\tau) &= \K_{3}(11/12,1)(\tau).
\end{split}
\end{equation}
The discussions above indicate that the $\K_{3}$ functions are only closed under the Kummer transformation \eqref{eq:k3kummer} in the three cases of \eqref{eq:k3tok3kmr}, unlike the Kummer transformation on the $\K_{2}$ functions.

We now consider the Atkin-Lehner involutions for the $\K_{2}, \K_{3},$ and $\K_{3}^{\kmr}$ functions in further detail.

\subsubsection{Atkin-Lehner Involutions}

An important Atkin-Lehner involution for the $\K_{2}$ functions is the map $W_{\K_{2}}: \tau \mapsto -1/\tau$ on $\Gamma(2)$, as mentioned earlier. In particular, the $\K_{2}$ functions satisfy 
\begin{equation}\label{eq:K2AL}
\K_{2}(r,s)(\tau)|_{W_{\K_{2}}} = 2^{4s-8r}i \cdot \tau^{3} \K_{2}(s-r,s)(\tau)
\end{equation}
for $(r,s) \in \mathbb{S}_{2}$, as discussed in \cite{hmm2}. Therefore, the $\K_{2}$ functions are closed under $W_{\K_{2}}$. 

\begin{remark}
The identities for Atkin-Lehner involutions in this section follow from the well-known transformation 
\begin{equation}\label{eq:etatransform}
\eta(-1/\tau) = \sqrt{\tau/i} \cdot \eta(\tau).
\end{equation}
\end{remark}

The important Atkin-Lehner involution for the $\K_{3}$ and $\K_{3}^{\kmr}$ functions is $W_{\K_{3}}: \tau \mapsto -1/(3 \tau)$ for $\Gamma_{0}(3)$. We now investigate to what extent the $\K_{3}$ and $\K_{3}^{\kmr}$ functions are invariant under $W_{\K_{3}}$, in analogy with the $\K_{2}$ case. Applying \eqref{eq:etatransform} to \eqref{eq:K3def} and \eqref{eq:K3kmr} gives
\begin{equation}\label{eq:K3ALtransform}
\K_{3}(r,1)(\tau)|_{W_{\K_{3}}} = 3^{\frac{9}{2}-6r}i \cdot \tau^{3}\K_{3}(1-r,1)(\tau)
\end{equation}
and 
\begin{equation}\label{eq:K3kummerALtransform}
\K_{3}^{\kmr}(r)(\tau)|_{W_{\K_{3}}} = 3^{\frac{1}{2}-6r}i \cdot \tau^{3} \eta(\tau)^{5+12r}\eta(3 \tau)^{1-12r}
\end{equation}
\begin{remark}
Note that $\eta(\tau)^{5+12r}\eta(3 \tau)^{1-12r}$ is not a  $\K_{3}^{\kmr}(r)(\tau)$ function due to the $r > 0$ requirement in \eqref{eq:S3kmrdef}.
\end{remark}

Therefore, the $\K_{3}$ functions are invariant under $W_{\K_{3}}$ by \eqref{eq:K3ALtransform} while the $\K_{3}^{\kmr}$ functions are not invariant under the same operator.

\bibliographystyle{plain}
\bibliography{ref2a}

@misc{ENRosenK2,
author={Rosen, Esme},
title={Transcendence of $_3{F}_2(1)$ Hypergeometric Series and \textit{L}-values of Modular Forms},
year={2024},
      eprint={2412.07054},
      archivePrefix={arXiv},
      primaryClass={math.NT},

}

@article{ENRosenK1,
author={Rosen, Esme},
title={Modular Forms and Certain $_2{F}_1(1)$ Hypergeometric Series},
year={2026},
journal={Proceedings of the American Mathematical Society}
}

@book {SerreAbelian,
    AUTHOR = {Serre, Jean-Pierre},
     TITLE = {Abelian {$l$}-adic representations and elliptic curves},
    SERIES = {Research Notes in Mathematics},
    VOLUME = {7},
      NOTE = {With the collaboration of Willem Kuyk and John Labute,
              Revised reprint of the 1968 original},
 PUBLISHER = {A K Peters, Ltd., Wellesley, MA},
      YEAR = {1998},
     PAGES = {199},
      ISBN = {1-56881-077-6},
   MRCLASS = {11G05 (11F33 11F80)},
  MRNUMBER = {1484415},
}

@incollection {Katz09,
    AUTHOR = {Katz, Nicholas M.},
     TITLE = {Another look at the {D}work family},
 BOOKTITLE = {Algebra, arithmetic, and geometry: in honor of {Y}u. {I}.
              {M}anin. {V}ol. {II}},
    SERIES = {Progr. Math.},
    VOLUME = {270},
     PAGES = {89--126},
 PUBLISHER = {Birkh\"{a}user Boston, Boston, MA},
      YEAR = {2009},
   MRCLASS = {14D10 (11L05 14C30 14D05 14G10)},
  MRNUMBER = {2641188},
MRREVIEWER = {Vasile Br\^{\i}nz\u{a}nescu},
       DOI = {10.1007/978-0-8176-4747-6\_4},
       URL = {https://doi-org.libezp.lib.lsu.edu/10.1007/978-0-8176-4747-6_4},
}

@article {HypMot,
    AUTHOR = {Roberts, David P. and Rodriguez Villegas, Fernando},
     TITLE = {Hypergeometric Motives},
   JOURNAL = {Notices Amer. Math. Soc.},
  FJOURNAL = {Notices of the American Mathematical Society},
    VOLUME = {69},
      YEAR = {2022},
    NUMBER = {6},
     PAGES = {914--929},
      ISSN = {0002-9920},
   MRCLASS = {01A70 (55-03 57-03)},
  MRNUMBER = {1691561},
}

@book {KatzESDE,
    AUTHOR = {Katz, Nicholas M.},
     TITLE = {Exponential sums and differential equations},
    SERIES = {Annals of Mathematics Studies},
    VOLUME = {124},
 PUBLISHER = {Princeton University Press, Princeton, NJ},
      YEAR = {1990},
     PAGES = {xii+430},
      ISBN = {0-691-08598-6; 0-691-08599-4},
   MRCLASS = {14D10 (11L03 11T23 14G15)},
  MRNUMBER = {1081536},
MRREVIEWER = {Hernando Enrique Sierra-Morales},
       DOI = {10.1515/9781400882434},
       URL = {https://doi-org.libezp.lib.lsu.edu/10.1515/9781400882434},
}

@article {lennon1,
    AUTHOR = {Lennon, Catherine},
     TITLE = {Trace formulas for {H}ecke operators, {G}aussian
              hypergeometric functions, and the modularity of a threefold},
   JOURNAL = {J. Number Theory},
  FJOURNAL = {Journal of Number Theory},
    VOLUME = {131},
      YEAR = {2011},
    NUMBER = {12},
     PAGES = {2320--2351},
      ISSN = {0022-314X},
   MRCLASS = {11F25 (11G20 11G40 11T24 33E50)},
  MRNUMBER = {2832827},
MRREVIEWER = {Jenny G. Fuselier},
       DOI = {10.1016/j.jnt.2011.05.005},
       URL = {https://doi-org.libezp.lib.lsu.edu/10.1016/j.jnt.2011.05.005},
}

@article {Allen,
    AUTHOR = {Allen, Michael},
     TITLE = {On some hypergeometric supercongruence conjectures of {L}ong},
   JOURNAL = {Ramanujan J.},
  FJOURNAL = {Ramanujan Journal. An International Journal Devoted to the
              Areas of Mathematics Influenced by Ramanujan},
    VOLUME = {61},
      YEAR = {2023},
    NUMBER = {3},
     PAGES = {957--987},
      ISSN = {1382-4090,1572-9303},
   MRCLASS = {11A07 (11B65 11F33 11S80 33C20)},
  MRNUMBER = {4599661},
       DOI = {10.1007/s11139-022-00665-2},
       URL = {https://doi.org/10.1007/s11139-022-00665-2},
}

@article {BBG-Ramanujan,
    AUTHOR = {Berndt, Bruce C. and Bhargava, S. and Garvan, Frank G.},
     TITLE = {Ramanujan's theories of elliptic functions to alternative
              bases},
   JOURNAL = {Trans. Amer. Math. Soc.},
  FJOURNAL = {Transactions of the American Mathematical Society},
    VOLUME = {347},
      YEAR = {1995},
    NUMBER = {11},
     PAGES = {4163--4244},
      ISSN = {0002-9947},
   MRCLASS = {33E05 (11F27 33C05 33D10)},
  MRNUMBER = {1311903},
MRREVIEWER = {J. Borwein},
       DOI = {10.2307/2155035},
       URL = {https://doi-org.libezp.lib.lsu.edu/10.2307/2155035},
}

@book {AAR,
    AUTHOR = {Andrews, George E. and Askey, Richard and Roy, Ranjan},
     TITLE = {Special functions},
    SERIES = {Encyclopedia of Mathematics and its Applications},
    VOLUME = {71},
 PUBLISHER = {Cambridge University Press, Cambridge},
      YEAR = {1999},
     PAGES = {xvi+664},
      ISBN = {0-521-62321-9; 0-521-78988-5},
   MRCLASS = {33-01 (33-02)},
  MRNUMBER = {1688958},
MRREVIEWER = {Bruce C. Berndt},
       DOI = {10.1017/CBO9781107325937},
       URL = {https://doi.org/10.1017/CBO9781107325937},
}

@article {Ahlgren-Ono-CalabiYau,
    AUTHOR = {Ahlgren, Scott and Ono, Ken},
     TITLE = {Modularity of a certain {C}alabi-{Y}au threefold},
   JOURNAL = {Monatsh. Math.},
  FJOURNAL = {Monatshefte f\"ur Mathematik},
    VOLUME = {129},
      YEAR = {2000},
    NUMBER = {3},
     PAGES = {177--190},
      ISSN = {0026-9255},
     CODEN = {MNMTA2},
   MRCLASS = {11G40 (11F20 11F23 11F80 14J32)},
  MRNUMBER = {1746757 (2001b:11059)},
MRREVIEWER = {Noriko Yui},
       DOI = {10.1007/s006050050069},
       URL = {http://dx.doi.org/10.1007/s006050050069},
}

@article {BCM,
    AUTHOR = {Beukers, Frits and Cohen, Henri and Mellit, Anton},
     TITLE = {Finite hypergeometric functions},
   JOURNAL = {Pure Appl. Math. Q.},
  FJOURNAL = {Pure and Applied Mathematics Quarterly},
    VOLUME = {11},
      YEAR = {2015},
    NUMBER = {4},
     PAGES = {559--589},
      ISSN = {1558-8599},
   MRCLASS = {11T24 (11L05 14G05 14M25 33C20 33C80)},
  MRNUMBER = {3613122},
MRREVIEWER = {Daniel Barsky},
       DOI = {10.4310/PAMQ.2015.v11.n4.a2},
       URL = {http://dx.doi.org/10.4310/PAMQ.2015.v11.n4.a2},
}

@book {BB,
    AUTHOR = {Borwein, Jonathan M. and Borwein, Peter B.},
     TITLE = {Pi and the {AGM}},
    SERIES = {Canadian Mathematical Society Series of Monographs and
              Advanced Texts},
    VOLUME = {4},
      NOTE = {A study in analytic number theory and computational complexity,
              Reprint of the 1987 original,
              A Wiley-Interscience Publication},
 PUBLISHER = {John Wiley \& Sons, Inc., New York},
      YEAR = {1998},
     PAGES = {xvi+414},
      ISBN = {0-471-31515-X},
   MRCLASS = {11Y60 (11B65 68Q25)},
  MRNUMBER = {1641658},
}

@article {BBG,
    AUTHOR = {Borwein, Jonathan M. and Borwein, Peter B. and Garvan, Frank. G.},
     TITLE = {Some cubic modular identities of {R}amanujan},
   JOURNAL = {Trans. Amer. Math. Soc.},
  FJOURNAL = {Transactions of the American Mathematical Society},
    VOLUME = {343},
      YEAR = {1994},
    NUMBER = {1},
     PAGES = {35--47},
      ISSN = {0002-9947},
   MRCLASS = {11B65 (11F27 33D10)},
  MRNUMBER = {1243610},
MRREVIEWER = {George E. Andrews},
       DOI = {10.2307/2154520},
       URL = {https://doi-org.libezp.lib.lsu.edu/10.2307/2154520},
}

@article {Win3X,
    AUTHOR = {Fuselier, Jenny and Long, Ling and Ramakrishna, Ravi and
              Swisher, Holly and Tu, Fang-Ting},
     TITLE = {Hypergeometric functions over finite fields},
   JOURNAL = {Mem. Amer. Math. Soc.},
  FJOURNAL = {Memoirs of the American Mathematical Society},
    VOLUME = {280},
      YEAR = {2022},
    NUMBER = {1382},
     PAGES = {},
      ISSN = {0065-9266},
      ISBN = {978-1-4704-5433-3; 978-1-4704-7282-5},
   MRCLASS = {11T23 (11-02 11F80 11S40 11T24 33)},
  MRNUMBER = {4493579},
       DOI = {10.1090/memo/1382},
       URL = {https://doi-org.libezp.lib.lsu.edu/10.1090/memo/1382},
}

@article {Greene,
    AUTHOR = {Greene, John},
     TITLE = {Hypergeometric functions over finite fields},
   JOURNAL = {Trans. Amer. Math. Soc.},
  FJOURNAL = {Transactions of the American Mathematical Society},
    VOLUME = {301},
      YEAR = {1987},
    NUMBER = {1},
     PAGES = {77--101},
      ISSN = {0002-9947},
     CODEN = {TAMTAM},
   MRCLASS = {11T21 (11L05 33A35)},
  MRNUMBER = {879564 (88e:11122)},
MRREVIEWER = {Ronald J. Evans},
       DOI = {10.2307/2000329},
       URL = {http://dx.doi.org/10.2307/2000329},
}

@article {LR,
    AUTHOR = {Long, Ling and Ramakrishna, Ravi},
     TITLE = {Some supercongruences occurring in truncated hypergeometric
              series},
   JOURNAL = {Adv. Math.},
  FJOURNAL = {Advances in Mathematics},
    VOLUME = {290},
      YEAR = {2016},
     PAGES = {773--808},
      ISSN = {0001-8708},
   MRCLASS = {33C20 (33E50)},
  MRNUMBER = {3451938},
MRREVIEWER = {Rupam Barman},
       DOI = {10.1016/j.aim.2015.11.043},
       URL = {https://doi.org/10.1016/j.aim.2015.11.043},
}

@article {Long,
    AUTHOR = {Long, Ling},
     TITLE = {Hypergeometric evaluation identities and supercongruences},
   JOURNAL = {Pacific J. Math.},
  FJOURNAL = {Pacific Journal of Mathematics},
    VOLUME = {249},
      YEAR = {2011},
    NUMBER = {2},
     PAGES = {405--418},
      ISSN = {0030-8730},
   MRCLASS = {33C20 (11B65)},
  MRNUMBER = {2782677},
MRREVIEWER = {Robert B. Osburn},
       DOI = {10.2140/pjm.2011.249.405},
       URL = {https://doi.org/10.2140/pjm.2011.249.405},
}

@incollection {K-Z,
    AUTHOR = {Kontsevich, Maxim and Zagier, Don},
     TITLE = {Periods},
 BOOKTITLE = {Mathematics unlimited---2001 and beyond},
     PAGES = {771--808},
 PUBLISHER = {Springer, Berlin},
      YEAR = {2001},
   MRCLASS = {11-02 (11F67 11G40 11G55)},
  MRNUMBER = {1852188},
MRREVIEWER = {F. Beukers},
}

@article {LTYZ,
    AUTHOR = {Long, Ling and Tu, Fang-Ting and Yui, Noriko and Zudilin,
              Wadim},
     TITLE = {Supercongruences for rigid hypergeometric {C}alabi-{Y}au
              threefolds},
   JOURNAL = {Adv. Math.},
  FJOURNAL = {Advances in Mathematics},
    VOLUME = {393},
      YEAR = {2021},
     PAGES = {Paper No. 108058, 49},
      ISSN = {0001-8708},
   MRCLASS = {11F33 (11T24 14G10 14J32 14J33 33C20)},
  MRNUMBER = {4330088},
MRREVIEWER = {Hidenori Katsurada},
       DOI = {10.1016/j.aim.2021.108058},
       URL = {https://doi-org.libezp.lib.lsu.edu/10.1016/j.aim.2021.108058},
}

@article {McCarthy-Papanikolas,
    AUTHOR = {McCarthy, Dermot and Papanikolas, Matthew A.},
     TITLE = {A finite field hypergeometric function associated to
              eigenvalues of a {S}iegel eigenform},
   JOURNAL = {Int. J. Number Theory},
  FJOURNAL = {International Journal of Number Theory},
    VOLUME = {11},
      YEAR = {2015},
    NUMBER = {8},
     PAGES = {2431--2450},
      ISSN = {1793-0421},
   MRCLASS = {11F46 (11F11 11G20 11T24 33E50)},
  MRNUMBER = {3420754},
MRREVIEWER = {Charles Helou},
       DOI = {10.1142/S1793042115501134},
       URL = {https://doi.org/10.1142/S1793042115501134},
}

@article {LPSX,
    AUTHOR = {Long, Ling and Plaza, Rafael and Sin, Peter and Xiang, Qing},
     TITLE = {Characterization of intersecting families of maximum size in
              {$PSL(2, q)$}},
   JOURNAL = {J. Combin. Theory Ser. A},
  FJOURNAL = {Journal of Combinatorial Theory. Series A},
    VOLUME = {157},
      YEAR = {2018},
     PAGES = {461--499},
      ISSN = {0097-3165},
   MRCLASS = {05D05},
  MRNUMBER = {3780424},
MRREVIEWER = {Norihide Tokushige},
       DOI = {10.1016/j.jcta.2018.03.006},
       URL = {https://doi.org/10.1016/j.jcta.2018.03.006},
}

@ARTICLE{LLT,
   author = {{Li}, Wen-Ching Winnie and {Long}, Ling and {Tu}, Fang-Ting},
    title = "{Computing special $L$-values of certain modular forms with complex multiplication}",
  journal = {SIGMA 14 (2018), 090},
 primaryClass = "math.NT",
     year = 2018,
    month = aug,
   adsurl = {http://adsabs.harvard.edu/abs/2018arXiv180306072L}
}

@article {LLT2,
    AUTHOR = {Li, Wen-Ching Winnie and Long, Ling and Tu, Fang-Ting},
     TITLE = {A {W}hipple {$_7F_6$} formula revisited},
   JOURNAL = {La Matematica},
  FJOURNAL = {La Matematica},
    VOLUME = {1},
      YEAR = {2022},
    NUMBER = {2},
     PAGES = {480--530},
   MRCLASS = {11F11 (11F67 11F80 33C20)},
  MRNUMBER = {4445932},
       DOI = {10.1007/s44007-021-00015-6},
       URL = {https://doi-org.libezp.lib.lsu.edu/10.1007/s44007-021-00015-6},
}

@article {Ono98,
    AUTHOR = {Ono, Ken},
     TITLE = {Values of {G}aussian hypergeometric series},
   JOURNAL = {Trans. Amer. Math. Soc.},
  FJOURNAL = {Transactions of the American Mathematical Society},
    VOLUME = {350},
      YEAR = {1998},
    NUMBER = {3},
     PAGES = {1205--1223},
      ISSN = {0002-9947},
   MRCLASS = {11T24 (33C99)},
  MRNUMBER = {1407498},
MRREVIEWER = {Ronald J. Evans},
       DOI = {10.1090/S0002-9947-98-01887-X},
       URL = {https://doi-org.libezp.lib.lsu.edu/10.1090/S0002-9947-98-01887-X},
}

@incollection {Salerno-Dwork,
    AUTHOR = {Salerno, Adriana},
     TITLE = {An algorithmic approach to the {D}work family},
 BOOKTITLE = {Women in numbers 2: research directions in number theory},
    SERIES = {Contemp. Math.},
    VOLUME = {606},
     PAGES = {83--100},
 PUBLISHER = {Amer. Math. Soc., Providence, RI},
      YEAR = {2013},
   MRCLASS = {11Y16 (14G10 33C20)},
  MRNUMBER = {3204293},
MRREVIEWER = {\'{I}sabel Pirsic},
       DOI = {10.1090/conm/606/12141},
       URL = {https://doi-org.libezp.lib.lsu.edu/10.1090/conm/606/12141},
}

@article {McCarthy-Dwork,
    AUTHOR = {McCarthy, Dermot},
     TITLE = {The number of {$\mathbb{F}_p$}-points on {D}work hypersurfaces
              and hypergeometric functions},
   JOURNAL = {Res. Math. Sci.},
  FJOURNAL = {Research in the Mathematical Sciences},
    VOLUME = {4},
      YEAR = {2017},
     PAGES = {Paper No. 4, 15},
      ISSN = {2522-0144},
   MRCLASS = {11G25 (11S80 11T24 33C99 33E50)},
  MRNUMBER = {3630724},
MRREVIEWER = {Rupam Barman},
       DOI = {10.1186/s40687-017-0096-y},
       URL = {https://doi-org.libezp.lib.lsu.edu/10.1186/s40687-017-0096-y},
}

@article {McCarthy,
    AUTHOR = {McCarthy, Dermot},
     TITLE = {Transformations of well-poised hypergeometric functions over finite fields},
   JOURNAL = {Finite Fields Appl.},
  FJOURNAL = {Finite Fields and their Applications},
    VOLUME = {18},
      YEAR = {2012},
    NUMBER = {6},
     PAGES = {1133--1147},
      ISSN = {1071-5797},
   MRCLASS = {11T24 (11L99 33C20)},
  MRNUMBER = {3019189},
MRREVIEWER = {Jenny G. Fuselier},
       DOI = {10.1016/j.ffa.2012.08.007},
       URL = {http://dx.doi.org/10.1016/j.ffa.2012.08.007},
}

@incollection {Zagier-top-diff,
    AUTHOR = {Zagier, Don},
     TITLE = {The arithmetic and topology of differential equations},
 BOOKTITLE = {European {C}ongress of {M}athematics},
     PAGES = {717--776},
 PUBLISHER = {Eur. Math. Soc., Z\"{u}rich},
      YEAR = {2018},
   MRCLASS = {11-02 (11Gxx 14H10 14J10 14J33)},
  MRNUMBER = {3890449},
}

@incollection {Dawsey-McCarthy,
    AUTHOR = {Dawsey, Madeline Locus and McCarthy, Dermot},
     TITLE = {Hypergeometric functions over finite fields and modular forms:
              a survey and new conjectures},
 BOOKTITLE = {From operator theory to orthogonal polynomials, combinatorics,
              and number theory---a volume in honor of {L}ance
              {L}ittlejohn's 70th birthday},
    SERIES = {Oper. Theory Adv. Appl.},
    VOLUME = {285},
     PAGES = {41--56},
 PUBLISHER = {Birkh\"{a}user/Springer, Cham},
      YEAR = {[2021] \copyright 2021},
      ISBN = {978-3-030-75424-2; 978-3-030-75425-9},
   MRCLASS = {33E50 (11F33)},
  MRNUMBER = {4367461},
       DOI = {10.1007/978-3-030-75425-9\_4},
       URL = {https://doi.org/10.1007/978-3-030-75425-9_4},
}

@article {Yu-Dwork,
    AUTHOR = {Yu, Jeng-Daw},
     TITLE = {Variation of the unit root along the {D}work family of
              {C}alabi-{Y}au varieties},
   JOURNAL = {Math. Ann.},
  FJOURNAL = {Mathematische Annalen},
    VOLUME = {343},
      YEAR = {2009},
    NUMBER = {1},
     PAGES = {53--78},
      ISSN = {0025-5831},
   MRCLASS = {14J32 (14F30 14G10)},
  MRNUMBER = {2448441},
MRREVIEWER = {Antonio Rojas Le\'{o}n},
       DOI = {10.1007/s00208-008-0265-9},
       URL = {https://doi-org.libezp.lib.lsu.edu/10.1007/s00208-008-0265-9},
}

@incollection {Zagier-modularform,
    AUTHOR = {Zagier, Don},
     TITLE = {Elliptic modular forms and their applications},
 BOOKTITLE = {The 1-2-3 of modular forms},
    SERIES = {Universitext},
     PAGES = {1--103},
 PUBLISHER = {Springer, Berlin},
      YEAR = {2008},
   MRCLASS = {11F11 (11-02 11E45 11F20 11F25 11F27 11F67)},
  MRNUMBER = {2409678},
MRREVIEWER = {Rainer Schulze-Pillot},
       DOI = {10.1007/978-3-540-74119-0_1},
       URL = {https://doi-org.libezp.lib.lsu.edu/10.1007/978-3-540-74119-0_1},
}

@book {Katz90,
    AUTHOR = {Katz, Nicholas M.},
     TITLE = {Exponential Sums and Differential Equations},
    SERIES = {Annals of Mathematics Studies},
    VOLUME = {124},
 PUBLISHER = {Princeton University Press, Princeton, NJ},
      YEAR = {1990},
     PAGES = {xii+430},
      ISBN = {0-691-08598-6; 0-691-08599-4},
   MRCLASS = {14D10 (11L03 11T23 14G15)},
  MRNUMBER = {1081536 (93a:14009)},
MRREVIEWER = {Hernando Enrique Sierra-Morales},
}

@incollection {Ramanujan-pi,
    AUTHOR = {Ramanujan, S.},
     TITLE = {Modular equations and approximations to {$\pi$} [{Q}uart. {J}.
              {M}ath. {\bf 45} (1914), 350--372]},
 BOOKTITLE = {Collected papers of {S}rinivasa {R}amanujan},
     PAGES = {23--39},
 PUBLISHER = {AMS Chelsea Publ., Providence, RI},
      YEAR = {2000},
      ISBN = {0-8218-2076-1},
   MRCLASS = {01A75},
  MRNUMBER = {2280849},
       URL = {},
}

@book {Ono-WebofModularity,
    AUTHOR = {Ono, Ken},
     TITLE = {The web of modularity: arithmetic of the coefficients of
              modular forms and {$q$}-series},
    SERIES = {CBMS Regional Conference Series in Mathematics},
    VOLUME = {102},
 PUBLISHER = {Conference Board of the Mathematical Sciences, Washington, DC;
              by the American Mathematical Society, Providence, RI},
      YEAR = {2004},
     PAGES = {viii+216},
      ISBN = {0-8218-3368-5},
   MRCLASS = {11F30 (05A30 11F11 11F37)},
  MRNUMBER = {2020489},
MRREVIEWER = {Kevin\ L.\ James},
}

@article{DMpaley,

     AUTHOR = {Dawsey, Madeline Locus and McCarthy, Dermot},
     TITLE = {Generalized {P}aley graphs and their complete subgraphs of orders three and four},
   JOURNAL = {Research in the Mathematical Sciences},
    VOLUME = {8},
      YEAR = {2021},
    NUMBER = {18},
}

@article{Mor3F2,
    author = {Eric Mortenson},
    title = {SUPERCONGRUENCES FOR TRUNCATED $_{n+1}\text{F}_{n}$
HYPERGEOMETRIC SERIES WITH APPLICATIONS
TO CERTAIN WEIGHT THREE NEWFORMS},
    journal = {Proc. Amer. Math. Soc.},
    volume = {13},
    year = {2004},
    number = {2},
}

@article{HMM1,
      AUTHOR = {Allen, Michael and Grove, Brian and Long, Ling and Fang-Ting Tu},
     TITLE = {The {E}xplicit {H}ypergeometric-{M}odularity {M}ethod {I}},
   JOURNAL = {Advances in Mathematics},
    VOLUME = {478},
      YEAR = {2025},
    NUMBER = {110411},
}

@article{HMM2,
      AUTHOR = {Allen, Michael and Grove, Brian and Long, Ling and Fang-Ting Tu},
     TITLE = {The {E}xplicit {H}ypergeometric-{M}odularity {M}ethod {II}},
   JOURNAL = {Research in the Mathematical Sciences},
    VOLUME = {12},
      YEAR = {2025},
    NUMBER = {84},
}

@article{originalPaley,

     AUTHOR = {R.E.A.C. Paley},
     TITLE = {On orthogonal matrices},
   JOURNAL = {Mathematics and Physics},
    VOLUME = {12},
     YEAR = {1933},
     PAGES = {311--320},
}

@article{evansPaley,

     AUTHOR = {R.J. Evans and J.R. Pulham and J. Sheehan},
     TITLE = {On the number of complete subgraphs contained in certain graphs},
   JOURNAL = {Combinatorial Theory Series B},
    VOLUME = {30},
     YEAR = {1981},
     NUMBER = {2},
     PAGES = {364-371},
}

@article{gPaley,

     AUTHOR = {T.K. Lim and C. Praeger},
     TITLE = {On generalized {P}aley graphs and their automorphism groups},
   JOURNAL = {Michigan Mathematics Journal},
    VOLUME = {58},
     YEAR = {2009},
     NUMBER = {1},
     PAGES = {293-308},
}

@article{ms1,
      TITLE={Orbits of {F}inite {F}ield {H}ypergeometric {F}unctions and {C}omplete {S}ubgraphs of {G}eneralized {P}aley {G}raphs}, 
      AUTHOR={Dermot McCarthy and Mason Springfield},
      JOURNAL = {Involve},
      VOLUME = {17},
      YEAR = {2024},
      NUMBER = {2},
      PAGES = {355-362}
}

@article{ms2,
      TITLE={Transitive subtournaments of $k$-th power {P}aley digraphs and improved lower bounds for {R}amsey numbers}, 
      AUTHOR={Dermot McCarthy and Mason Springfield},
      JOURNAL = {Graphs and Combinatorics},
      VOLUME = {40},
      YEAR = {2024},
      NUMBER = {71},
}

@article{bg1,

     AUTHOR = {Anwita Bhowmik and Rupam Barman},
     TITLE = {On a {P}aley-{T}ype Graph on $\mathbb{Z}_{n}$},
   JOURNAL = {Graphs and Combinatorics},
    VOLUME = {38},
     YEAR = {2022},
     NUMBER = {4},
}

@article{bg2,

     AUTHOR = {Anwita Bhowmik and Rupam Barman},
     TITLE = {Number of complete subgraphs of {P}eisert graphs and finite field hypergeometric functions},
   JOURNAL = {Research in Number Theory},
    VOLUME = {10},
     YEAR = {2024},
     NUMBER = {26},
}

@article{bg3,
      title = {Cliques of orders three and four in the {P}aley-type graphs}, 
      author = {Anwita Bhowmik and Rupam Barman},
      year = {2024},
      journal = {Graphs and Combinatorics},
}

@article{bg4,

     AUTHOR = {Anwita Bhowmik and Rupam Barman},
     TITLE = {Hypergeometric {F}unctions for {D}irichlet {C}haracters and {P}eisert-{L}ike {G}raphs on $\mathbb{Z}_{n}$},
   JOURNAL = {La Mathematica},
    VOLUME = {2},
     YEAR = {2023},
     PAGES = {992-1021},
}

@article{koike2dim,
title = {Orthogonal matrices obtained from hypergeometric series over finite fields and elliptic curves over finite fields},
author={Masao Koike},
journal = {Hiroshima Mathematics Journal},
volume = {25},
pages = {43-52},
year = {1995},
}

@article{momentint,
title = {{M}oments of elliptic integrals and critical {$L$}-values},
author={Rogers, M. and Wan, J.G. and Zucker, I.J.},
journal = {The Ramanujan Journal},
volume = {37},
pages = {113-130},
year = {2015},
}

@article{special3vals,
title = {{T}he special values of {$L$}-functions at $s=1$ of theta products of weight 3},
author={Ito, Ryojun},
journal = {Research in Number Theory},
volume = {5},
number = {30},
year = {2019},
}

@article{Wiles,
  title={Modular Elliptic Curves and {F}ermat's {L}ast {T}heorem},
  author={Wiles, Andrew},
  journal={Annals of Mathematics},
  volume={141},
  number={3},
  pages={443-551},
  year={1995},
}

@article{TaylorWiles,
  title={Ring-theoretic properties of certain {H}ecke algebras},
  author={Richard Taylor and Andew Wiles},
  journal={Annals of Mathematics},
  volume={141},
  number={3},
  pages={553-572},
  year={1995},
}

@article{wild3adic,
  title={On the modularity of elliptic curves over $\mathbb{Q}$: {W}ild $3$-adic exercises},
  author={Christophe Breuil and Brian Conrad and Fred Diamond and Richard Taylor},
  journal={Journal of the American Mathematical Society},
  volume={14},
  number={4},
  pages={843-939},
  year={2001},
}

@article{conditionalST1,
 author = {Michael Harris and Nick Shepherd-Barron and Richard Taylor},
 journal = {Annals of Mathematics},
 number = {2},
 pages = {779--813},
 title = {A family of {C}alabi--{Y}au varieties and potential automorphy},
 volume = {171},
 year = {2010}
}

@article{conditionalST2,
 author = {Richard Taylor},
 journal = {Publications mathématiques},
 pages = {183-239},
 title = {Automorphy for some $l$-adic lifts of automorphic mod $l$ {G}alois representations. {II}},
 volume = {108},
 year = {2008}
}

@article{conditionalST3,
 author = {Laurent Clozel and Michael Harris and Richard Taylor},
 journal = {Publications mathématiques},
 pages = {1-181},
 title = {Automorphy for some $l$-adic lifts of automorphic mod $l$ {G}alois representations},
 volume = {108},
 year = {2008}
}

@article{unconditionalST,
 author = {Tom Barnet-Lamb and David Geraghty and Michael Harris and Richard Taylor},
 journal = {Publications of the Research Institute for Mathematical Sciences},
 publisher = {Kyoto University},
 pages = {29-98},
 title = {A Family of {C}alabi–-{Y}au Varieties and Potential Automorphy {II}},
 number = {1},
 volume = {47},
 year = {2011}
}

@article{GoodsonDwork,
title = {A complete hypergeometric point count formula for {D}work hypersurfaces},
journal = {Journal of Number Theory},
volume = {179},
pages = {142-171},
year = {2017},
author = {Heidi Goodson},
}

@article{BarmanDwork, 
title={Counting points on {D}work hypersurface and $p$-adic hypergeometric function}, 
volume={94}, 
number={2}, 
journal={Bulletin of the Australian Mathematical Society}, 
author={Barman, Rupam and Rahman, Hasanur and Saikia, Neelam}, 
year={2016}, 
pages={208–216}}

@article{Borweincubic,
 author = {J. M. Borwein and P. B. Borwein},
 journal = {Transactions of the American Mathematical Society},
 number = {2},
 pages = {691--701},
 publisher = {American Mathematical Society},
 title = {A {C}ubic {C}ounterpart of {J}acobi's {I}dentity and the {AGM}},
 volume = {323},
 year = {1991}
}

@incollection{jonesPaley,
author="Jones, Gareth A.",
title="Paley and the {P}aley {G}raphs",
booktitle="Isomorphisms, Symmetry and Computations in Algebraic Graph Theory",
year="2020",
publisher="Springer International Publishing",
pages="155--183",
}

\end{document}